\documentclass[11pt]{article}
\usepackage{graphicx,,psfrag}
\usepackage{delarray}
\usepackage{amsmath, latexsym, amsfonts, amssymb, amsthm, amscd}
\usepackage{graphicx}
\usepackage{graphicx,,psfrag}
\usepackage{amssymb,amsmath,amsthm,enumerate}
\usepackage{color}
\newcommand{\Tr}{\operatorname{Tr}}
\newcommand{\1}{1\!\!{\sf I}}

 \newtheorem{theorem}{Theorem}
 \newtheorem{proposition}{Proposition}
\newtheorem{remark}{Remark}

\title{Spectrum of deformed random matrices and free probability}
\author{M. Capitaine\thanks{CNRS, Institut de Math\'ematiques de Toulouse, Universit\'e Paul Sabatier, 31062 Toulouse, France}, C. Donati-Martin\thanks{Laboratoire de Math\'ematiques de Versailles, UVSQ, CNRS, Universit\'e Paris-Saclay 78035 Versailles, France}}
\date{}
\begin{document}
\maketitle
\begin{abstract}The aim of this paper is to show how free probability theory sheds light on spectral properties of deformed matricial models and provides a unified understanding of various asymptotic  phenomena such as  spectral measure description, localization and fluctuations of extremal eigenvalues, eigenvectors behaviour.
\end{abstract}

\section{Introduction}

The sum or product of large independent random matrices are models of interest in applications such as mathematical finance, wireless communication, statistical learning
(see for example Johnstone \cite{Jo}).   Indeed, in these contexts, one matrix is seen as the signal with significant parameters and the other one as    the noise. The general question  is to know whether the observation 
of signal plus noise can give access to the significant parameters.

  Pastur \cite{Pastur} and Marchenko-Pastur \cite{MP} investigated  the limiting spectral distributions of respectively additive deformation of Wigner matrices  or multiplicative deformation  of Wishart matrices and by the way provide a description of the global asymptotic behavior of the spectrum. 
Concerning the behavior of largest  eigenvalues,
the pionner works goes back to Furedi-Komlos \cite{FurKom81} where it is shown that for a non centered   Wigner matrix, the largest eigenvalue separates from  the bulk of eigenvalues and has Gaussian fluctuations.
In 2005, Baik, Ben Arous and P\'ech\'e \cite{BBP} proved a phase transition phenomenon for the behavior of the largest eigenvalue of  the so-called non white Wishart matrices (that is  multiplicative  perturbations of a Wishart matrix). The multiplicative perturbation 
they considered is a finite rank perturbation of the identity matrix. The authors pointed out a threshold (depending on the covariance of the sample vectors) giving two different behaviors~: either the largest eigenvalue sticks to the bulk and fluctuates according to the Tracy Widom distribution, either the largest eigenvalue separates from the bulk and has Gaussian fluctuations.
The phenomenon described  is  called the BBP transition.
Later, a phase transition phenomenon was  proved for  more general finite rank  deformations of     both matrices of iid type or  unitarily invariant models \cite{BaikSil06, Peche, BN1, BN2, LV} and for general spiked sample covariance matrices \cite{BY}. Some papers investigated eigenvectors of finite rank  deformations \cite{Paul, BN1,BN2}. \\

 It is well known since Voiculescu's works \cite{Voiculescu91} that free probability sheds light on the global asymptotic spectrum of  the sum or product of large independent random matrices. 
In the last five years, several articles have highlighted the role of the free analytic subordination functions (introduced by  Voiculescu \cite{voic-fish1, Voiculescu00}, Biane \cite{Biane98}) 
in the analysis of both global and local behaviour of eigenvalues and eigenspaces of polynomials in asymptotically free random matrices. Thus, the eigenvalue distribution of any polynomial in asymptotically free random matrices can be found via its Cauchy-Stieltjes transform through a direct computation of the subordination functions associated to a free additive convolution of matrix-valued free random variables. We refer to \cite{BSTV} and \cite{BMS} for the analysis of selfadjoint polynomials, and to \cite{BSS} for the determination of Brown measure for non-selfadjoint polynomials. More surprisingly, the analytic subordination property for free convolutions turned out to be fundamental in the study of outliers of spiked deformed models (eigenvalues of the deformed matrix converging outside the support of the limiting spectral distribution) and 
  allows the study to be performed on arbitrary classical  deformed  models, while explaining the first order asymptotic in the BBP phenomenon \cite{CDFF, BBCF,MCJTP,MC2014, Shly}. The same subordination functions provide (via their derivatives) the limiting behaviour of the eigenvectors associated to the outliers (\cite{BN1, BN2, MCJTP, BBCF}).
Very recently, the local laws of sums of independent unitarily invariant random matrices is studied in (\cite{BES}) using subordination.
\\

In this paper, we investigate three classical deformed  models. It turns out that  various properties of the spectrum, i.e. convergence (or fluctuations)  of eigenvalues and eigenvectors,  of these models can be presented in an unified way, using the framework of free probability, free convolutions and free subordination functions.    Therefore, this paper will focus on free probability and  take another look of various results in the literature from this point of view without being exhaustive.


In the rest of this introduction, we describe the three models of deformed random matrices: additive deformation of a Wigner or unitarily invariant  matrix, multiplicative deformation of a Wishart type  or unitarily invariant matrix and information plus noise type model. We then present the results on the limiting spectral distributions and the seminal works on the convergence of extreme eigenvalues, for  finite rank  perturbations. 

Section 2 is devoted to free probability : we first provide some background on the theory and then focus on free convolutions and subordinations properties which will be central in the study of general deformed models.
In Section 3, we give a free probabilistic interpretation of the limiting spectral distribution (LSD) of Section 1 and describe the support of the LSD, in terms of subordinations functions.
Section 4 is devoted to the behavior of outliers  of general spikes models. The location of these outliers as well as the  norm of the projection of corresponding eigenvectors onto the spikes eigengenspaces of the perturbation   are described by means of  the subordinations functions.
In the last section, we describe some results on fluctuations of outliers and eigenvalues at soft edges of the LSD. The first results was obtained for Gaussian models with finite rank deformations.
For general deformations, one can get universal result by considering mobile edges of the support of the  Girko's type deterministic equivalent measure obtained by replacing the LSD of $A_N$ by its empirical  spectral measure  in the free convolution describing the  limiting LSD of the deformed model.  
We also provide a free probabilistic interpretation of the appearance of different rates or different asymptotic distributions  for fluctuations as uncovered by 
Johansson \cite{J} and Lee and Schnelli \cite{LeeSchnelli}.

\vspace{.3cm}
\noindent
{\bf Notations:} \\
For a Hermitian matrix $H$ of size $N \times N$, we denote by $\{\lambda_1(H) \geq \lambda_2(H) \geq \ldots \geq \lambda_N(H) \}$ the set of eigenvalues of $H$, ranked in decreasing order, and by 
$\displaystyle \mu_H = \frac{1}{N} \sum_{i=1}^N \delta_{\lambda_i(H)}$ its
 spectral measure. The Stieltjes transform of $\mu_H$ is defined for $z \in \mathbb C \backslash \mathbb R$ by 
$$ g_{\mu_H} (z) = \int_{\mathbb R} \frac{1}{z-x} \mathrm{d}\mu_H(x) = \frac{1}{N} \Tr(zI-H)^{-1}.$$
If $\mu_H$ converges weakly to some probability measure $\mathbb{P}$ as $N$ goes to infinity, $\mathbb{P}$ is called the limiting spectral measure (LSD) of $H$ and $ \mbox{supp}(\mathbb{P})$ denotes its support.

\subsection{Classical Random matricial models and deformations}
We first review classical models in random matrix theory and the pionner works on the behavior of their spectrum.
\subsubsection{Wigner matrices}
Wigner matrices are real symmetric or complex Hermitian random matrices whose entries are independent (up to the symmetry condition).
They were introduced by Wigner in the fifties, in connection with nuclear physics.
In this paper, we will consider Hermitian Wigner matrices of the following form :
$$W_N= \frac{1}{\sqrt{N}} X_N$$
 where $(X_N)_{ii}$, $\sqrt{2} \Re e((X_N)_{ij})_{i<j}$, $\sqrt{2}
\Im m ((X_N)_{ij})_{i<j}$
are  i.i.d,
with distribution $\tau$ with  variance $\sigma^2$ and mean zero.\\
  If $\tau = {\cal N}(0,\sigma^2)$, $ W_N=:W_N^G$ is a G.U.E.-matrix.
 The first result concerns the behavior of the spectral measure.
  \begin{theorem}[\cite{Wigner55,Wigner58}]
 $$\mu_{W_N}:= \frac{1}{N} \sum_{i=1}^N \delta_{\lambda_i(W_N)} \rightarrow \mu_{sc}  \mbox{~~a.s when~~} N \rightarrow + \infty  $$
 where 
$$ \frac{d \mu_{sc}}{dx}(x)=  \frac{1}{2 \pi \sigma^2}
\sqrt{4\sigma^2 - x^2} \, 1_{[-2\sigma,2 \sigma]}(x)$$
\end{theorem}
The second result establishes the behavior of the extremal eigenvalues.
\begin{theorem} [\cite{BaiYin88}] If $\int x^4 \mathrm{d}\tau(x)< + \infty$, then  $$\lambda_{1}(W_N) \rightarrow 2 \sigma
\mbox{~and~} \lambda_{N}(W_N) \rightarrow -2 \sigma
\mbox{~a.s when~} N \rightarrow + \infty.$$
\end{theorem}
\subsubsection{Sample covariance matrices}
Random matrices first appeared in mathematical statistics in the 1930s with the works of Hsu, Wishart and others.
They considered sample covariance matrices of the form:
\begin{equation}\label{wishart}S_N= \frac{1}{p} X_N  X_N^*\end{equation}
where $X_N$ is a random matrix with i.i.d. entries. \\
We shall assume in the following (complex case) that $N \leq p(N) $,
 $X_N$ is a $N\times p(N)$ matrix, $\Re e((X_N)_{ij}), \Im  m ((X_N)_{ij}), \, i=1,\ldots,N$, $j=1, \ldots, p$ are  i.i.d,
with distribution $\tau$ with variance $\frac{1}{2}$ and mean zero. Note that the spectra of $ \frac{1}{p} X_N  X_N^*$ and $\frac{1}{p} X_N^* X_N$ differ by $\vert p-N\vert$ zero eigenvalues. If $\tau$ is Gaussian, $S_N=:S_N^G$ is  a L.U.E. matrix.\\

The behavior of the spectral measure for large size ( $p=p(N)$ tends to $\infty$ as  $N$ tends to $\infty$)  was handled in the seminal work of Marchenko-Pastur.
\begin{theorem}[\cite{MP}]
If $c_N:=\frac{N}{p(N)} \rightarrow c\in ]0;1]$ when N $\rightarrow \infty$, 
$$\mu_{S_N} \rightarrow \mu_{\operatorname{MP}}  \mbox{~~a.s when~~} N \rightarrow + \infty  $$
where 
 $$\mu_{\operatorname{MP}}(dx)=  \frac{1}{2 \pi c x}
\sqrt{(b-x)(x-a)} \, 1_{[a,b]}(x) \mathrm{d}x$$
where $a= (1-{\sqrt{c}})^2$, $b= (1+{\sqrt{c}})^2 .$ 
\end{theorem} 
For extremal eigenvalues of sample covariance matrices, we have : 
\begin{theorem}[\cite{Ge,BaiSilYi,YBK}]
 If $\int x^4 \mathrm{d}\tau(x)< + \infty$,
$$\lambda_{1}(S_N) \rightarrow  (1+{\sqrt{c}})^2
\mbox{~a.s when~} N \rightarrow + \infty,$$
$$\lambda_{N}(S_N) \rightarrow  (1-{\sqrt{c}})^2
\mbox{~a.s when~} N \rightarrow + \infty.$$
\end{theorem}

\subsubsection{Deformations\label{siid}}
Motivated by statistics, wireless communications or imaging one may consider deformations of these models.

\noindent
Let $A_N$ be a deterministic  matrix. One may wonder how the spectrum of a classical random model is impacted by the following perturbations.
\begin{itemize}
\item[i)] Additive perturbation of a Wigner matrix :   $A_N$ is a  $N \times N$  Hermitian matrix and $W_N$ a Wigner matrix (see section 1.1.1), 
$$M_N = W_N +A_N.$$
\item[ii)]Multiplicative perturbation:  $A_N $ is  non negative Hermitian $N \times N$  matrix and $S_N$ a sample covariance matrix defined as in  section 1.1.2,
$$M_N={A_N}^{\frac{1}{2}}S_N {A_N}^{\frac{1}{2}}.$$
\item[iii)] Information-Plus-Noise type model:  $A_N$, $X_N$  are  rectangular $N\times p$  matrices  with $X_N$ a random matrix with i.i.d. entries defined as in  section 1.1.2, $\sigma$ is some positive real number.
$$M_N= (\sigma \frac{ X_N}{\sqrt{p}}+{A_N})(\sigma\frac{X_N}{\sqrt{p}}+{A_N})^*.$$
\end{itemize}
In the litterature, these three kinds of deformations have been also considered for isotropic models, replacing in i) and ii) $W_N$, $S_N$ by $UBU^*$ with $U$ Haar distributed and $B$ deterministic, and in iii) $X_N$ by a random matrix whose distribution is biunitarily invariant.
\subsection{Convergence of spectral measures}
This section is devoted to the study of the limiting spectral measure of the deformed matrix $M_N$ for the different models i) to iii). 
\subsubsection{Finite rank deformations}
When $A_N$ is a finite rank deformation of the null matrix in case i) and iii) (resp. of the identity matrix in case ii)),  
  the limiting spectral measure is not affected by the deformation.
 \begin{proposition} We assume that  ${\rm rank}(A_N) = r$, $r$ fixed, independent of $N$ in case i) and iii), (resp. ${\rm rank}(A_N-I) =r$ in case ii)). Then, 
when $N$ goes to infinity,
\begin{itemize}
\item In case i), $\mu_{M_N} \rightarrow \mu_{sc}$,
\item In case ii) and iii), $\mu_{M_N} \rightarrow \mu_{\operatorname{MP}}$.
\end{itemize}
\end{proposition}
\noindent
This follows from the rank inequalities (see \cite[Appendix A.6]{BS10}) where for a matrix $A$, $F^A$ denotes the cumulative distribution function of the spectral measure $\mu_A$ and the norm is the supremum norm: \\
- Let $A$ and $B$ two $N \times N$ Hermitian matrices. Then,
$$\Vert F^A - F^B \Vert \leq \frac{1}{N} \rm{rank}(A-B),$$
- Let $A$ and $B$ two $N\times p$ matrices. Then,
$$\Vert F^{AA^*} - F^{BB^*} \Vert \leq \frac{1}{N} \rm{rank}(A-B).$$
.

\subsubsection{General deformations}
We now assume that for some probability $\nu$, the spectrum of the deformation $A_N$ satisfies 
$$ \begin{array}{ll}  \mu_{A_N} \underset{N \to + \infty}{\longrightarrow} \nu & \mbox{~~for cases i) and ii)}\\
\mu_{A_N A_N^*} \underset{N \to + \infty}{\longrightarrow} \nu & \mbox{~~for case iii)}. \end{array}$$ 
The following theorem is a review of the pionner results concerning the limiting spectral distribution, in the three  deformed models. The spectral distribution is characterized via an  equation satisfied by its Stieltjes transform. In Section \ref{freeinterpretation}, we will give a free probabilistic interpretation of these distributions.
\begin{theorem}\label{convmes} Convergence of the spectral measure $\mu_{M_N}$ 
\begin{itemize}
\item[i)] Deformed Wigner matrices (\cite{Pastur}, \cite{AGZ10} Theorem 5.4.5) \\
$$\mu_{M_N} \underset{N \to + \infty}{\longrightarrow} \mu_{1} \; {\rm weakly} $$
with
\begin{equation} \label{TSwigner} \forall z \in \mathbb{C}^+,~~g_{\mu_{1}}(z)=\int \frac{1}{z-\sigma^2 g_{\mu_{1}}(z)-t}\mathrm{d}\nu(t). \end{equation}

\item[ii)] Sample covariance matrices 
(\cite{MP, S95})
$$\mu_{M_N} \underset{N \to + \infty}{\longrightarrow} \mu_{2} \; {\rm  weakly} $$
with 
\begin{equation} \label{TScovariance} \forall z \in \mathbb{C}^+,~~g_{\mu_{2}}(z)=\int \frac{1}{z-t(1-c +czg_{\mu_{2}}(z))}\mathrm{d}\nu(t). \end{equation}
\item[iii)] Information-Plus-Noise type  matrices (\cite{DozierSilver})
$$ \mu_{M_N}\underset{N \to + \infty}{\longrightarrow} \mu_{3} \; {\rm weakly}  $$
with
\begin{equation} \label{TSnoise}\hspace*{-1,3cm} \forall z \in \mathbb{C}^+,~~g_{\mu_{3}}(z)=\int \frac{1}{(1- c \sigma^2g_{\mu_{3}}(z))z- \frac{ t}{1-  c\sigma^2g_{\mu_{3}}(z)} - \sigma^2 (1-c)}\mathrm{d}\nu(t). \end{equation}
\end{itemize}
\end{theorem}
\begin{remark}\label{remarque1}
The limiting measures $\mu_1,\mu_2,\mu_3$ are deterministic. Note that they are not always explicit. They are  {\bf universal} (do not depend on the distribution of the entries of $X_N$) and only depend on $A_N$ through the limiting spectral measure $\nu$.
\end{remark}
\begin{remark}
If $\nu = \delta_0$ in \eqref{TSwigner}, we recover the equation satisfied by the Stieltjes transform of the semicircular distribution. The same is true for the Marchenko-Pastur distribution in \eqref{TScovariance} with $\nu = \delta_1$ and in \eqref{TSnoise} with $\nu = \delta_0$.
\end{remark}
\begin{remark}
Such functional equations for Stieltjes transforms have been obtained for deformations of unitarily invariant models by \cite{PV} and \cite{V}.
\end{remark}
\subsection{Convergence of extreme eigenvalues} \label{seminal}
We now present the seminal works on the behavior of the largest (or smallest) eigenvalues  of classical models with finite rank  deformation.
As  we have seen above, the limiting behavior of the spectral measure is not modified by a deformation $A_N$ of finite rank (or such that $I_N- A_N$ is of finite rank in the multiplicative case). This is no longer true for the extremal eigenvalues. \\
The following results deal with finite rank perturbations and  Gaussian type or unitarily invariant models : we give the precise statement for the convergence of the largest eigenvalue in the Gaussian case and refer to the original papers in the unitarily invariant case.

\subsubsection{Multiplicative deformations}
A complete study (convergence, fluctuations) of the behavior of the largest eigenvalue of a deformation of a Gaussian sample covariance matrix was considered in a paper of Baik-Ben Arous and P\'ech\'e  where they exhibit a striking  phase transition phenomenon for the largest eigenvalue, according to the value of the spiked eigenvalues of the deformation. They considered the following sample covariance matrix :
$${M_N= A_N^{1/2}S_N^G A_N^{1/2}}$$
where $S_N^G$ is a L.U.E. matrix as defined in Section 1.1.2 and the perturbation $A_N$\footnote{In the Gaussian case, we may assume, without loss of generality, that the perturbation is diagonal}  is given by 
$${A_N= \mbox{diag~} (\underbrace{1, \ldots,1}_{N-r \mbox{~times}}, \pi_1,\ldots, \pi_r)}$$
where $r$ is fixed, independent of $N$ and 
$\pi_1 \geq \ldots \geq \pi_r >0$ are  fixed, independent of $N$, such that  for all $ i \in \{1,\ldots,r\}, ~ \pi_i \neq 1$. The $\pi_i$'s are called the spikes of $A_N$.
 
\begin{theorem}(BBP phase transition)\cite{BBP,BaikSil06}

  Let $\omega_c= 1 + \sqrt{c},$

 \begin{itemize} \item {If $\pi_1> \omega_c$}, $ \mbox{~a.s when~} N \rightarrow + \infty$ $$\lambda_1 \left(M_N\right) 
 \rightarrow  \pi_1 \left( 1 + \frac{c}{(\pi_1-1)}\right) > (1+{\sqrt{c}})^2.$$
Therefore the largest eigenvalue of $M_N$  is an {``outlier"} since it converges  outside the support of
the limiting empirical spectral distribution and then does not stick to the bulk.
 \item If $\pi_1 \leq  \omega_c$, $\mbox{~a.s when~} N \rightarrow + \infty$  $$\lambda_1 \left(M_N\right) 
 \rightarrow(1+{\sqrt{c}})^2. $$
 \end{itemize}
\end{theorem}
The same phenomenon of phase transition was established by Benaych-Georges and Nadakuditi \cite{BN1} in the case of a deformation of an unitarily invariant model of the form :
\begin{equation} \label{uniinv}
M_N=A_N^{1/2}U_N B_N U_N^*A_N^{1/2},
\end{equation}
where
\begin{itemize}
 \item $B_N$ is a deterministic $N\times N$ Hermitian non negative definite matrix such that:
\item $\mu _{B_N} := \frac{1}{N} \sum_{i=1}^N \delta _{\lambda _i(B_N)}$ 
weakly converges to some  probability measure $\mu $ with compact support $[a,b]$.
\item the smallest and largest eigenvalues of $B_N$ converge  to $a$ and $b$.
\item $U_N$ is a random $N\times N$ unitary matrix distributed according to Haar measure. 
\item $A_N - I_N$ is a deterministic Hermitian  matrix 
having $r$ non-zero eigenvalues $\theta_1 \geq \ldots \geq \theta_r > -1$
$r$, $\theta_i, i=1,\ldots r,$ fixed independent of $N$ .
\end{itemize}
\subsubsection{Additive deformations}
Let us consider $M_N=W_N^G+A_N$ where $W_N^G$ is a G.U.E. matrix as defined in Section 1.1.1 and $A_N$ is defined by 
$$A_N= \mbox{diag~} (\underbrace{0, \ldots,0}_{N-r \mbox{~times}},\theta_1, \ldots,\theta_r)$$
for some fixed  $r$, independent of $N$, and some fixed  $\theta_1\geq \cdots \geq  \theta_r$,  independent of $N$.

\noindent An analog of the BBP phase transition phenomenon for this model was obtained by P\'ech\'e \cite{Peche}.
 \begin{theorem}

\begin{itemize}

\item If $\theta_1 \leq \sigma$, then $\lambda_1(M_N) \underset{N \to + \infty}{\longrightarrow} 2 \sigma$ a.s..

\item If $\theta_1 >\sigma$, then $\lambda_1(M_N)\underset{N \to + \infty}{\longrightarrow}\rho_{\theta_1}$ a.s. 
with $\rho_{\theta_1}:= \theta_1+ \frac{\sigma^2}{\theta_1} > 2 \sigma$.
\end{itemize}
\end{theorem}
\noindent
A similar result has been obtained by Benaych-Georges and Nadakuditi \cite{BN1} for
$$M_N=U_N B_N U_N^*+A_N,$$
where $U_N$, $B_N$ satisfy the same hypothesis as in \eqref{uniinv} and $A_N$ is a finite rank perturbation of the null matrix.

\subsubsection{Information plus noise type matrices}
Such a phase transition phenomenon was established by Loubaton and Vallet in \cite{LV} for
 the singular values of a finite rank deformation of a Ginibre ensemble. \\
Let $X_N$ be a $N \times p$ rectangular matrix as defined in Section 1.1.2 with iid complex Gaussian entries, and $A_N$ be a finite rank perturbation of the null matrix with non zero eigenvalues $\theta_i$.
\begin{theorem}(\cite{LV})
Let $M_N=( \sigma \frac{ X_N}{\sqrt{p}}+A_N)(\sigma \frac{ X_N}{\sqrt{p}}+A_N)^*$.
Then, as $N \rightarrow + \infty$ and $N/p \rightarrow c\in ]0;1]$  , almost surely, 
$$\lambda_i(M_N)  ~~\longrightarrow \left\{ \begin{array}{ll} \frac{(\sigma^2 +\theta_i)(\sigma^2 c +\theta_i)}{\theta_i} \mbox{~~if~~} \theta_i > \sigma^2 \sqrt{c},\\ \sigma^2 (1+\sqrt{c})^2 \mbox{~~otherwise}. \end{array}\right.$$
\end{theorem}

Again, this result was extended in \cite{BN2} to the case $M_N = (V_N +A_N)(V_N + A_N)^*$ where $V_N$ is a biunitarily invariant matrix such that the empirical  spectral measure of $V_N V_N^*$ converges to a deterministic compactly supported measure $\mu$ with convergence of the largest (resp. smallest) eigenvalue of $V_NV_N^*$ to the right (resp. left) end of the support of $\mu$.
\subsection{Eigenvectors associated to outliers}
Now, one can wonder in the spiked deformed models, when some eigenvalues separate from the bulk, how the corresponding eigenvectors of the deformed model project onto those of the perturbation.
There are some pionneering results concerning finite
rank perturbations: \cite{Paul} in the real Gaussian sample covariance matrix setting,
and \cite{BN1,BN2} dealing with finite rank additive or multiplicative perturbations
of unitarily invariant matrices. Here is  the result for eigenvector projection corresponding to the largest outlier  for finite rank additive perturbation.
\begin{theorem}
Let $M_N=U_N B_N U_N^*+A_N$ where  $U_N$ is a Haar unitary matrix, $B_N$ satisfies the same hypothesis as in (\ref{uniinv}) and  $A_N$ has all but finitely many non zero eigenvalues $\theta_1>\ldots>\theta_J$.
Then, if  $\theta_1>1/\lim_{z\downarrow b} g_\mu( z)$, almost surely, 
$$\lambda_1(M_N) \underset{N \to + \infty}{\longrightarrow}  g_{\mu}^{-1}(1/\theta_1):=\rho_{\theta_1}$$
and for any $i =1, \ldots,J$, if $\xi$ is a  unit eigenvector associated to $\lambda_1(M_N)$,
$$ \Vert P_{\operatorname{Ker} (\theta_i I-A)} \xi \Vert^2 \underset{N \to + \infty}{\longrightarrow} -\frac{\delta_{i1}}{\theta_1^2 g_\mu^{'}(\rho_{\theta_1})}.$$
\end{theorem}

It turns out that the results of the above Section may be interpreted in terms of free probability theory and 
this very analysis allows  to extend the
results of Subsections 1.3 and 1.4 to non-finite rank deformations. Therefore, in the next section, we start by recalling some   background in free probability theory for the  reader's convenience.

\section{Free probability theory}
We refer to \cite{VDN92} for an
introduction to free probability theory.  \\

A non-commutative probability space is a unital algebra ${\cal A}$ over
$\mathbb{C}$, endowed with a linear functional $\phi : A \rightarrow
\mathbb{C}$ such
that $ \phi(1)=1$. Elements of ${\cal A}$ are called non-commutative  random variables.\\

If $(a_i)_{i=1,\ldots,q}$ is a family of non-commutative random variables in $({\cal A},\phi)$, 
the distribution $\mu_{(a_i)_{i=1,\ldots,q}}$ of $(a_i)_{i=1,\ldots,q}$ is the linear functional on the algebra
$\mathbb{C} \langle X_i \vert i=1,\ldots,q \rangle $ of polynomials in the non commutating variables $(X_i)_{i=1,\ldots,q}$
given by
$$\mu_{(a_i)_{i=1,\ldots,q}}(P) = \phi \left( P\left(\mu_{(a_i)_{i=1,\ldots,q}}\right) \right).$$
If ${\cal A}$ is a $C^*$-algebra endowed with a state $\phi$, then for any selfadjoint element $a$ in  ${\cal A}$,  there exists a measure $\nu_a$ on $\mathbb{R}$ such that,   for every polynomial P, we have
$$\mu_a(P)=\int P(t) \mathrm{d}\nu_a(t).$$
Then  we identify $\mu_a$ and $\nu_a$.\\

A family of unital subalgebras $\left({\cal A}_i\right)_{i=1, \ldots q}$ in $({\cal A},\phi)$ is freely independent
if for every $p \geq 1$, for every $(a_1,\ldots,a_p)$ such that, for every $k$ in $\{1,\ldots,p\}$, 
$\phi(a_k)=0$ and $a_k$ is in ${\cal A}_{i(k)}$ for some $i(k)$ in $\{1,\ldots,q\}$ with
$i(k)\neq i(k+1)$, then $\phi(a_1,\ldots,a_p)=0$.
Random variables are free in  $({\cal A},\phi)$ if the subalgebras they generate with 1
are freely independent.\\

 For each $n$ in ${\mathbb N}\setminus \{0\}$, let $(a_i^n)_{i=1,\ldots,q}$ be a family of noncommutative
random variables in a noncommutative probability space $({\cal A}_n,\phi_n)$. The
sequence of joint distributions $\mu_{(a_i^n)_{i=1,\ldots,q}}$ converges as $n$  tends to $+\infty$, if there exists
a distribution $\mu$ such that $\mu_{(a_i^n)_{i=1,\ldots,q}}(P)$ converges to $\mu(P)$ as $n$ tends to $+\infty$ for
every $P$ in $\mathbb{C} \langle X_i \vert i=1,\ldots,q \rangle $. $\mu$ is called the limit distribution of $(a_i^n)_{i=1,\ldots,q}$. If
$(a_i)_{i=1,\ldots,q}$ is a family of noncommutative random variables with distribution $\mu$,
we also say that $(a_i^n)_{i=1,\ldots,q}$ converges towards $(a_i)_{i=1,\ldots,q}$.
A family of noncommutative random variables $(a_i^n)_{i=1,\ldots,q}$ is said to be asymptotically
free as $n$ tends to $\infty$ if it has a limit distribution $\mu$ and if $(X_1,\ldots,X_q)$ are
free in $\left(\mathbb{C} \langle X_i \vert i=1,\ldots,q \rangle, \mu \right)$.\\

Additive and multiplicative free convolutions arise as natural analogues of classical convolutions
in the context of free probability theory. 
For two 
Borel probability measures $\mu$ and $\nu$ on the real line, one defines the free 
additive convolution $\mu\boxplus\nu$ as the distribution of $a+b$, where $a$ and $b$
are free self-adjoint random variables with distributions $\mu$ and
$\nu$, respectively. Similarly, if both $\mu,\nu$ are supported on $[0,+\infty)$, their free multiplicative convolution $\mu\boxtimes\nu$ is the distribution of
the product $ab$, where, as before, $a$ and $b$
are free positive  random variables with distributions 
$\mu$ and $\nu$, respectively. The product $ab$ of two
free positive random variables is usually not positive, but it has the same moments as the 
positive random variables $a^{1/2}ba^{1/2}$ and $b^{1/2}ab^{1/2}.$ We refer  to \cite{Voiculescu86,V2, Maassen, BercoviciVoiculescu} for the
definitions and main properties of free convolutions. In the following  sections,
we briefly recall the analytic approach developed in \cite{Voiculescu86,V2}
to calculate the additive and multiplicative free convolutions of compactly supported measures and  the analytical definition of the rectangular free convolution introduced by F. Benaych-Georges in \cite{FB}.   Finally, we present the fundamental analytic subordination
properties \cite{voic-fish1,Biane98,Voiculescu00, BBG} of these three convolutions.

\subsection{Free  convolution}

\subsubsection{Free  additive convolution}
The Stieltjes transform of a compactly supported probability
measure $\mu$ is conformal in the neighborhood of $\infty$, and
its functional inverse $g_{\mu}^{-1}$ is meromorphic at zero with
principal part $1/z$. The $R$-transform \cite{Voiculescu86} of $\mu$ is the convergent
power series defined by 
\[
R_{\mu}(z)=g_{\mu}^{-1}(z)-\frac{1}{z}.
\]
 The free additive convolution of two compactly supported probability
measures $\mu$ and $\nu$ is another compactly supported probability
measure characterized by the identity 
\[
R_{\mu\boxplus\nu}=R_{\mu}+R_{\nu}
\]
satisfied by these convergent power series.

\subsubsection{Multiplicative free convolution on $[0,+\infty)$}
Recall  the definition of the moment-generating function of a 
Borel probability measure $\mu$ on $[0,+\infty)$:
$$
\psi_\mu(z)=\int_{[0,+\infty)}\frac{zt}{1-zt}\,\mathrm{d}\mu(t), \quad z\in\mathbb C\setminus
\left\{z\in\mathbb C\colon\frac1z\in\text{supp}(\mu)\right\}.
$$
This function is related to the Cauchy-Stieltjes transform of $\mu$ via the relation
$$\psi_\mu(z)=\frac1zg_\mu\left(\frac1z\right)-1.$$ 
Recall also the so-called eta transform
$$\eta_\mu(z)=\frac{\psi_\mu(z)}{1+\psi_\mu(z)}.$$

The 
$\Sigma$-transform \cite{V2,BercoviciVoiculescu} of a compactly supported Borel probability measure  $\mu\ne\delta_0$ 
is the convergent power series defined by
$$
\Sigma_\mu(z)=\frac{\eta_\mu^{-1}(z)}{z},
$$
where $\eta_{\mu}^{-1}$ is the inverse of $\eta_\mu$ relative to composition.
The free multiplicative convolution of two compactly supported probability
measures $\mu\ne\delta_0\ne\nu$ is another compactly supported probability
measure characterized by the identity 
\[
\Sigma_{\mu\boxtimes\nu}(z)=\Sigma_{\mu}(z)\Sigma_{\nu}(z)
\]
in a neighbourhood of $0$.
\subsubsection{Rectangular free convolution}
Let $c$ be in $ ]0;1]$.
Let
$\tau$ be a probability measure on $\mathbb{R^+}$. Define for $z$ in $\mathbb{C}\setminus [0;+\infty[$, 
$$M_\tau (z)=\int_{\mathbb{R^+}} \frac{t^2z}{1-t^2 z} \mathrm{d}\tau(t), ~~{H_\tau^{(c)}(z):= z \left( cM_\tau (z) +1 \right) \left( M_\tau (z) +1 \right)},$$
$$\mbox{and~~} T^{(c)}(z) =(cz +1)(z+1).$$
The  transform $C_\tau$ \cite{FB} defined as follows is called the rectangular R-transform: $$C_\tau^{(c)} (z) = {T^{(c)}}^{-1} \left( \frac{z}{{H_\tau^{(c)}}^{-1}(z)} \right), ~~\mbox{for}~~ z  \mbox{~~small enough}.$$

The rectangular free convolution with ratio $c$  of two  probability measures $\mu$ and $\nu$ on $\mathbb{R^+}$
is the unique  probability measure   on $\mathbb{R^+}$ whose rectangular R-transform is the sum of the rectangular R-transforms of $\mu$ and $\nu$, and it is denoted by $\mu \boxplus_c \nu$. Then, we have for $z$ small enough,
$$C_{\mu \boxplus_c \nu}(z)=C_{\mu }(z)+C_{ \nu}(z).$$
\subsection{Analytic subordinations}

The analytic subordination phenomenon for free convolutions was first noted by Voiculescu in \cite{voic-fish1} for free 
additive convolution of compactly supported probability measures.  Biane \cite{Biane98}
extended the result to free additive convolutions of arbitrary probability measures on $\mathbb R$,
and also found a subordination result for multiplicative free convolution.
A new proof was
given later, using a fixed point theorem for analytic self-maps of the upper
half-plane \cite{BelBer07}. Note that such a subordination property  allows to give
a new definition of free additive convolution \cite{CG08a}. Finally, S. Belinschi, F. Benaych-Georges, and A. Guionnet \cite{BBG} established such a phenomenon for the rectangular free convolution.
\subsubsection{Free additive 
subordination property} Let us define  the reciprocal Cauchy-Stieltjes
transform $J_{\mu}(z)={1}/{g_{\mu}(z)}$, which is an analytic
self-map of the upper half-plane.  Given Borel probability measures $\mu$ and $\nu$ on 
$\mathbb R$, there exist two unique analytic functions $\omega_1,\omega_2\colon
\mathbb C^+\to\mathbb C^+$ such that 
\begin{enumerate}
\item $\lim_{y\to+\infty}\omega_j(iy)/iy=1$, $j=1,2$;
\item
\begin{equation}\label{EqSubord+}
\omega_1(z)+\omega_2(z)-z=J_\mu(\omega_1(z))=J_\nu(\omega_2(z))=J_{\mu\boxplus
\nu}(z),\quad z\in\mathbb C^+.
\end{equation}
\item In particular (see \cite{BelBer07}), for any $z\in\mathbb C^+\cup\mathbb R$ so that 
$\omega_1$ is analytic at $z$, $\omega_1(z)$ is the attracting fixed point of the self-map
of $\mathbb C^+$ defined by 
$$
 w\mapsto J_\nu(J_\mu(w)-w+z)-(J_\mu(w)-w).
$$
A similar statement, with $\mu,\nu$ interchanged, holds for $\omega_2$.
\end{enumerate}
In particular, according to (\ref{EqSubord+}), we have for any $z \in \mathbb C^+,$
\begin{equation}\label{subord1}
g_{\mu\boxplus\nu}(z)=g_\mu(\omega_1(z))=g_\nu(\omega_2(z)).
\end{equation}

\subsubsection{Multiplicative subordination property\label{subsubsec:eqX}} Given Borel probability measures $\mu,
\nu$ on $[0,+\infty)$, there exist two unique analytic functions $F_1,F_2
\colon\mathbb C\setminus[0,+\infty)\to\mathbb C\setminus[0,+\infty)$ so that 
\begin{enumerate}
\item $\pi>\arg F_j(z)\ge\arg z$ for $z\in\mathbb C^+$ and $j=1,2$;
\item
\begin{equation}\label{EqSubordX+}
\frac{F_1(z)F_2(z)}{z}=\eta_\mu(F_1(z))=\eta_\nu(\omega_2(z))=\eta_{\mu\boxtimes
\nu}(z),\quad z\in\mathbb C\setminus[0,+\infty).
\end{equation}
\item In particular (see \cite{BelBer07}), for any $z\in\mathbb C^+\cup\mathbb R$ so that 
$F_1$ is analytic at $z$, the point $h_1(z):=F_1(z)/z$ is the attracting fixed point of 
the self-map of ${\mathbb C}\setminus[0,+\infty)$ defined by
$$
 w\mapsto\frac{w}{\eta_\mu(zw)} 
\eta_\nu\left(\frac{\eta_\mu(zw)}{w}\right).
$$
A similar statement, with $\mu,\nu$ interchanged, holds for $\omega_2$.
\end{enumerate}
In particular (\ref{EqSubordX+}) yields
\begin{equation}\label{subord1'}
\psi_{\mu\boxtimes\nu}(z)=\psi_\mu(F_1(z))=\psi_\nu(F_2(z)).
\end{equation}

\begin{subsubsection}{Rectangular free subordination property } 

Let $c $ be in $ ]0;1]$.
Let $\mu $ and $\nu$ be two probability measures on $\mathbb{R^+}$.
Assume that the rectangular R-transform $C_\mu^{(c)}$ of $\mu$ extends analytically to $\mathbb{C}\setminus \mathbb{R}^+$; this happens for example if $\mu$ is $\boxplus_c$ infinitely divisible. 
 Then there exist two unique meromorphic functions $\Omega_1$, $\Omega_2$ on $\mathbb{C}\setminus \mathbb{R}^+$
so that $$H^{(c)}_\mu(\Omega_1(z))=H^{(c)}_\nu(\Omega_2(z)) =H^{(c)}_{\mu \boxplus_c \nu}(z),$$
$\Omega_j(\overline{z})=\overline{ \Omega_j(z)}$ and $\lim_{x \uparrow 0} \Omega_j(x)=0$, $j\in\{1;2\}$.

\end{subsubsection}
~~

The fonctions $\omega_i$, $F_i$ and $\Omega_i$ in the three previous subsections are called subordination functions.
\subsection{Asymptotic freeness of independent random matrices}
Free probability theory and random matrix theory are closely related.
Indeed the purely algebraic concept of free relation of
noncommutative random variables can be also modeled by random matrix ensembles
if the matrix size goes to infinity. Let  $(\Omega, {\cal F}, P)$ be a classical probability space and for every $N \geq 1$,
let us denote by  ${\cal M}_N$ the algebra of complex $N\times N$ matrices. Let 
${\cal A}_N$ be the algebra of $N\times N$ random matrices $(\Omega, {\cal F}, P) \rightarrow {\cal M}_N$.
 Define
$$\Phi_N : \left\{ \begin{array}{ll} {\cal A}_N \rightarrow
\mathbb{C} \\
A \rightarrow \frac{1}{N}\operatorname{Tr} A
\end{array}\right.$$
\noindent For each  $N \geq 1$,
$\left( {\cal A}_N , \Phi_N
\right)$
 is a non-commutative probability space and we will
consider random matrices in this non-commutative context.\\
 In his pioneering
work \cite{Voiculescu91}, Voiculescu shows  that independent Gaussian Wigner matrices converge in distribution
as their size goes to infinity to free semi-circular variables. The following result from Corollary 5.4.11 \cite{AGZ10} extends Voiculescu's seminal
observation.

\begin{theorem}\label{free}
 Let $\{D_N(i)
\}_{1\leq i \leq p}$ be a sequence of uniformly bounded real diagonal
matrices with empirical measure of diagonal elements converging to $\mu_i$,
$i = 1, \ldots , p$ respectively. Let $\{U_N
(i )\}_{1\leq i\leq p}$ be independent unitary matrices following
the Haar measure, independent from $\{D_N(i)
\}_{1\leq i \leq p}$.
\begin{itemize}
\item  The noncommutative variables $\{U_N(i) D_N(i) U_N(i)^*\}_{1\leq i\leq p}$ in the noncommutative
probability space $\left({\cal A}_N, \Phi_N \right)$ are  almost surely asymptotically
free, the law of the marginals being given by the $\mu_i$'s.
\item  The spectral distribution of $D_N(1)+U_N(2)D_N(2)U_N(2)^*$ converges weakly almost surely
to $\mu_1\boxplus\mu_2$ goes to infinity.
\item  Assume that $D_N(1)$ and $D_N(2)$
are nonnegative. Then, the empirical spectral measure of $(D_N(1))^{\frac{1}{2}}U_N(2)D_N(2)U_N(2)^*(D_N(1))^{\frac{1}{2}}$
converges weakly almost surely to $\mu_1 \boxtimes\mu_2$ as $N$ goes to infinity.
\end{itemize}
\end{theorem}
Thus, if $\mu$ is the eigenvalue distribution of a large selfadjoint random matrix $A$ and $\nu$ is the eigenvalue distribution of a large selfadjoint random matrix $B$ then, when $A$ and $B$ are in generic position,  $\mu\boxplus \nu$ is nearly the eigenvalue distribution 
of $A+B$. Similarly, dealing with nonnegative matrices  $\mu\boxtimes \nu$ is nearly the eigenvalue distribution 
of $A^{\frac{1}{2}}BA^{\frac{1}{2}} $.\\

{Similarly, for   independent rectangular $n\times N$ random matrices $A$ and $B$  such that $n/N \rightarrow c \in]0;1]$, when $A$ and $B$ are in generic position, F. Benaych-Georges \cite{FB} proved that rectangular free convolution with ratio $c$ provides a good understanding of the asymptotic global behaviour of 
 the singular values of
  $A+B$.
\begin{theorem}\label{rectfree}
Let $A$ and $B$ be  independent rectangular $N\times p$ random matrices such that $A$ or $B$
is invariant, in law, under multiplication, on the right and on the left, by any unitary matrix.
Assume that 
  there exists two laws $\mu$ and $\nu$
such that, for the weak convergence in probability, we
have  $\displaystyle{\frac{1}{N} \sum_{s~ 
\mbox{\small sing. val. of~}
A_N} \delta_s \rightarrow \mu}$,
~ $\displaystyle{\frac{1}{N} \sum_{s
\mbox{\small ~ sing. val. of ~}B_N} \delta_s \rightarrow \nu}$
as $n$ and $N$ goes to infinity with $N/p \rightarrow c \in]0;1]$.
Then 
$$\frac{1}{N} \sum_{s 
\mbox{\small ~sing. val. of~}
A_N+B_N} \delta_s \underset{N \to + \infty}{\longrightarrow} \mu \boxplus_c \nu,$$
 for the weak convergence in probability.
\end{theorem}

\section{Study of LSD of deformed ensembles through free probability theory } \label{freeinterpretation}
In this section, we take an other look of the LSD described in Section 1.2, in the light of free probability theory introduced in the previous Section. We also characterize the limiting supports in terms of free subordination functions. This will prove to be fundamental to understand the outliers phenomenon for spiked models in Section 4. We finish with an analysis of the different behaviors of the density at edges of the support of  free additive, multiplicative, rectangular convolutions  with semi-circular, Marchenko-Pastur and the square-root of Marchenko-Pastur distributions respectively. This provides the  rate for fluctuations of eigenvalues at edges as it will discussed in Section 5.
\subsection{Free probabilistic interpretation of LSD}
As noticed in Remark \ref{remarque1}, the limiting spectral distributions of the deformed models investigated in Section \ref{siid} are  {\bf universal} in the sense that they do not depend on the distribution of the entries of $X_N$. Therefore, choosing Gaussian entries and applying Theorem \ref{free} and Theorem \ref{rectfree}, we readily get the following free probabilistic interpretation of the limiting measures as well as of the equations satisfied by the limiting Stieltjes transforms in Theorem \ref{convmes}.
\begin{itemize}
\item[] $\bullet$ {Deformed Wigner matrices} 
$$\mu_{M_N} \underset{N \to + \infty}{\longrightarrow} {\mu_{1}} \mbox{ weakly}, ~~
{\mu_{1}=\mu_{sc}\boxplus{\nu}}$$
\item[] $\bullet$ {Sample covariance matrices}
$$\mu_{M_N} \underset{N \to + \infty}{\longrightarrow} {\mu_{2}} \mbox{ weakly},~~
{\mu_{2}=\mu_{\operatorname{MP}}\boxtimes{\nu}}$$
\item[]  $\bullet$ {Information-Plus-Noise type  matrices}
 $$\mu_{M_N} \underset{N \to + \infty}{\longrightarrow} {\mu_{3}} \mbox{ weakly},~~
{\mu_{3}=(\sqrt{\mu_{\operatorname{MP}}}\boxplus_c \sqrt{\nu})^2}.$$
\end{itemize}

The equations (\ref{TSwigner}), (\ref{TScovariance}) and (\ref{TSnoise}) satisfied by the limiting Stieltjes transforms  correspond to free subordination properties and exhibit the subordination functions $\omega_{\mu_{sc},\nu}$ with respect to the semi-circular distribution $\mu_{sc}$ for the free additive convolution, $F_{\mu_{\operatorname{MP}},\nu}$ with respect to the Marchenko-Pastur distribution $\mu_{\operatorname{MP}}$ for the free multiplicative  convolution, and $\Omega_{\mu_{\operatorname{MP}},\nu}$
with respect to the pushforward of the Marchenko-Pastur distribution by the square root function $\sqrt{\mu_{\operatorname{MP}}}$ for the rectangular free  convolution.\\

 $\bullet$ {Deformed Wigner matrices} 
{\begin{equation} \forall z \in \mathbb{C}^+,~~ g_{\mu_{1}}(z)=\int \frac{1}{z-\sigma^2g_{\mu_{1}}(z)-t}\mathrm{d}\nu(t)={g_\nu}({\omega_{\mu_{sc},\nu}}(z)) .\nonumber \end{equation}}
\vspace*{-0.2cm} $$\omega_{\mu_{sc},\nu}(z)=z-\sigma^2 g_{\mu_{1}}(z).$$

$\bullet$ {Sample covariance matrices}
{\begin{equation} \forall z \in \mathbb{C}^+,~~g_{\mu_{2}}(z)=\int \frac{1}{z-t(1-c +czg_{\mu_{2}}(z))}\mathrm{d}\nu(t).\nonumber \end{equation}}
  $$\rightarrow  ~~~~\psi_{\mu_{1}}\left(\frac{1}{z}\right)= {\psi_{\nu}}({F_{\mu_{\operatorname{MP}},\nu}}\left(\frac{1}{z}\right))$$
$${\psi_\tau (z)}=\int \frac{tz}{1-tz}d\tau(t)={\frac{1}{z}g_\tau(\frac{1}{z}) -1},$$ \vspace*{-0.2cm}$$
{F_{\mu_{\operatorname{MP}},\nu}(z)= z-cz +cg_{ \mu_{1}}(\frac{1}{z}).}$$

 $\bullet$ {Information-Plus-Noise type  matrices}
$${\mu_{3}=(\sqrt{\mu_{\operatorname{MP}}}\boxtimes_c \sqrt{\nu})^2}$$
{\begin{equation} \hspace*{-1,3cm} \forall z \in \mathbb{C}^+,~~g_{\mu_{3}}(z)=\int \frac{1}{(1- c \sigma^2g_{\mu_3}(z))z- \frac{ t}{1-  c\sigma^2g_{\mu_{3}}(z)} - \sigma^2(1-c)}\mathrm{d}\nu(t).\nonumber \end{equation}}
$$ \rightarrow  ~~{H^{(c)}_{\sqrt{\mu_{3}}}}\left(\frac{1}{z}\right) ={H^{(c)}_{\sqrt{\nu}}}\left({\Omega_{\mu_{\operatorname{MP}},\nu}}\left(\frac{1}{z}\right)\right) $$
$${H^{(c)}_{\sqrt{\tau}} (z) =\frac{c}{z} g_\tau(\frac{1}{z})^2 +(1-c) g_\tau(\frac{1}{z})},$$
$${\Omega_{\mu_{\operatorname{MP}},\nu}(z)=\frac{1}{
\frac{1}{z} (1-  c\sigma^2g_{ \mu_{3}}(\frac{1}{z}))^2 - (1-c)\sigma^2 (1- c \sigma^2g_{\mu_{3}}(\frac{1}{z}))}}.$$

\subsection{Limiting supports of LSD} \label{sectionsupport}
For each deformed model introduced in Section \ref{siid} involving i.i.d entries, several authors studied the limiting support \cite{Biane97b, CHOISILVER, DozierSilver2, VLM, LV, MC2014}.

It turns out that in each case there is a one to one correspondance involving the subordination  functions between  the complement of the support of the limiting spectral measure   and some set  in the complement of the limiting support of the deformation, as follows.\\

$\bullet$ \underline{Deformed Wigner} $\mu_1= \mu_{sc} \boxplus \nu$\\

 $~~~~~~~~~~~\mathbb{R}\setminus  {\rm supp}(\mu_1)
 \begin{array}{cc}\stackrel{{\varphi_1}}{\longrightarrow} \\
 \stackrel{\longleftarrow}{{\phi_1}}
 \end{array}  ~{{\cal O}_1}\subset \mathbb{R}\setminus  {\rm supp}(\nu) ,$
$${\cal O}_1:=\{ u \in \mathbb{R} \setminus {\rm supp} (\nu), \phi'_1(u)>0\}$$
 $$u \in \mathbb{R} \setminus {\rm supp} (\nu),~~{\phi_1(u~)=u+\sigma^2 g_{ \nu}(u)}.$$
$$x \in \mathbb{R} \setminus {\rm supp} (\mu_1), ~~{\varphi_1(x~)=x-\sigma^2 g_{ \mu_1}(x)}.$$
~~

$\bullet$ \underline{Sample covariance matrices } $\mu_2=\mu_{\operatorname{MP}} \boxtimes \nu$
\begin{equation}\label{caracsuppsample} {\mathbb{R}\setminus \{{\rm supp} (\mu_2)\}}
 \begin{array}{cc}\stackrel{{\varphi_2}}{\longrightarrow} \\
 \stackrel{\longleftarrow}{{\phi_2}}
 \end{array}  ~{{\cal O}_2}\subset \mathbb{R}\setminus  \{{\rm supp} (\nu) \} ,\end{equation}
$$ {\cal O}_2=\left\{ u \in  ^c \{\mbox{supp}( \nu)\},\, \phi'_2(u)>0\right\} $$

  $$u\in \mathbb{R}\setminus  {\rm supp} (\nu),~~{\phi_2(u)=u + cu \int \frac{t}{u-t}\mathrm{d}\nu(t).}$$ 
$$x \in \mathbb{R} \setminus {\rm supp} (\mu_2),~~\varphi_2(x)=~\left\{\begin{array}{ll} \frac{x}{(1-c) +cx g_{\mu_2}(x)} \mbox{~if ~} c<1\\
\frac{1}{ g_{\mu_2}(x)} \mbox{~~~~~~~~~~~if ~} c=1. \end{array}
\right.$$
Note that  $\varphi_2$ is well defined on $ \mathbb{R} \setminus {\rm supp} (\mu_2)$ since its denominator never vanishes according to Lemma 6.1 in \cite{BS10}.\\

$\bullet$ \underline{Information-Plus-Noise type model} $\mu_3=(\sqrt{\mu_{\operatorname{MP}}}\boxtimes_c \sqrt{\nu})^2$\\

$\mathbb{R}\setminus \rm{supp} (\mu_3)
 \begin{array}{cc}\stackrel{{ \varphi_3}}{\longrightarrow} \\
 \stackrel{\longleftarrow}{{\phi_3}}
 \end{array}  ~{{\cal O}_3}\subset \mathbb{R}\setminus  {\rm supp} (\nu) ,$
$$
 {\cal O}_3=\left\{ u \in  \mathbb{R}\setminus \mbox{supp}(\nu),  \phi_{3}^{'}(u) >0, g_\nu(u) >-\frac{1}{\sigma^2c}\right\}.$$
\begin{equation}\label{phi3} u \in \mathbb{R} \setminus {\rm supp} (\nu),~~{\phi_3(u)=u(1+c\sigma^2g_\nu(u))^2 + \sigma^2 (1-c)(1+c\sigma^2g_\nu(u)) }\end{equation}
$$x \in \mathbb{R} \setminus {\rm supp} (\mu_3),~~\varphi_3(x)=x (1-  c\sigma^2g_{ \mu_{3}}(x))^2 - (1-c)\sigma^2 (1- c \sigma^2g_{\mu_{3}}(x))$$
~~

Note that $\varphi_1$ corresponds to the extension of  $\omega_{\mu_{sc},\nu}$ on $\mathbb{R}\setminus {\rm supp} (\mu_1)$ and  for $i=2,3$, $\varphi_i$ coincides on 
$ \mathbb{R} \setminus \{{\rm supp} (\mu_i) \cup\{0\}\}$, with the extension of $z\mapsto 1/F_{\mu_{\operatorname{MP}},\nu}(1/z)$ and 
$z\mapsto 1/\Omega_{\mu_{\operatorname{MP}},\nu}(1/z)$ respectively. \\

The above characterization of the support are explicitely given in \cite{Biane97b,MC2014} for $\mu_1$ and $\mu_3$. 
Now, it can be deduced for $\mu_2$ by the following arguments. In a $W^*$-probability space endowed with a faithful state, the support of the distribution of a random variable $x$ corresponds to the spectrum of $x$. Thus, 
considering $\mu_2$ as the distribution of $b^{1/2} a b^{1/2}$ where  $a$ and $b$ are free bounded operators whose distributions are $\mu_{\operatorname{MP}}$ and $\nu$ respectively, one can easily see that for $c<1$, $0$ belongs to the support of $\mu_{\operatorname{MP}}\boxtimes \nu$ if and only if $0$ belongs to the support of $\nu$. The latter equivalence and 
 Lemma 6.1 in \cite{BS10} readily yield (\ref{caracsuppsample}).\\

When the support of $\nu$ has a finite number of connected components, we
have the following description of the support of the $\mu_i$'s in terms of a finite union
of closed disjoint intervals.
\begin{theorem}\label{descsupp} \cite{CDFF, MC2014}\label{compconnexes}
 Assume that the support of $\nu$ is a finite union of disjoint (possibly degenerate) closed bounded intervals. For any $i=1,3$,
there exists a nonnul integer number  $p$ and  $u_1< v_1<u_2<\ldots <u_p < v_p$ (depending on $i$) such that $${\cal O}_{i}=]-\infty,u_1[ \;\cup_{l=1}^{p-1} \; ]v_l,u_{l+1}[\; \cup \;
]v_p,+\infty[.$$
We have $$\mbox{supp}(\nu) \subset \cup_{l=1}^{p} [u_l,v_l]$$
and for  each $l \in \{1,\ldots,p\}$,  $[u_l,v_l]\cap\mbox{supp}(\nu) \neq \emptyset$.

\noindent Moreover,
$$
 \mbox{supp}(\mu_{i})
=\cup_{l=1}^p [\phi_{i}(u_l^-),\phi_{i}(v_l^+)],$$
with 
$$\phi_{i}(u_1^-) < \phi_{i}(v_1^+)< \phi_{i}(u_2^-) < \phi_{i}(v_2^+)
<\cdots< \phi_{i}(u_p^-) < \phi_{i}(v_p^+),$$
where $ \phi_{i}(u_l^-) =\lim_{u\uparrow u_l} \phi_{i}(u)$ and $ \phi_{i}(v_l^+) =\lim_{u\downarrow v_l} \phi_{i}(u)$.

Finally, for  each $l \in \{1,\ldots,p\}$, \begin{equation}\label{palier}\mu _{i }([\phi _{i }(u_l^-), \phi_{i }(v_l^+)]) = \nu ([u_l, v_l]).\end{equation}

\end{theorem}
Using the characterization of the support (\ref{caracsuppsample}}), Remark 3.6 in \cite{MCJTP} and the fact that  from \cite{B} the only possible mass of $\mu_2= \mu_{\operatorname{MP}}\boxtimes \nu$ is at zero,
one may check that the above  result still holds for $\mu_2$ allowing $u_1=v_1=0$ or $\phi_2(u_1)=\phi_2(v_1)=0$ in Theorem \ref{descsupp}. Note that the latter cases
occur only when $\nu$ has a Dirac mass at zero since from \cite{B}, $\mu_{\operatorname{MP}}\boxtimes \nu(\{0\})=\max (\mu_{\operatorname{MP}}(\{0\}), \nu(\{0\}))$ and therefore, since $c \leq 1$, $\mu_2$ has a Dirac mass at zero if and only if $\nu$ has a Dirac mass at zero. (\ref{palier}) can be seen as a  consequence of the matricial exact separation phenomenon described in Section 4.1.1 b) below letting $N$ go to infinity.

\subsubsection{Behavior of the density at edges}
P. Biane  proved in \cite{Biane97b} that $\mu_1=\mu_{sc}\boxplus\nu$ has a continuous density. Choi and Silverstein \cite{CHOISILVER} and Dozier and Silverstein \cite{DozierSilver2} proved respectively  that, away from zero, $\mu_2=\mu_{\operatorname{MP}} \boxtimes \nu$ and   $\mu_3=(\sqrt{\mu_{\operatorname{MP}}}\boxtimes_c \sqrt{\nu})^2$ possess a continuous density. Let us denote any of theses densities by $p$.\\
Using the notations of Theorem \ref{compconnexes},  we have $$\mbox{supp}(\nu) \subset \cup_{l=1}^{p} [u_l;v_l]=\mathbb{R}\setminus  {\cal O}_{i}$$
and for  each $l \in \{1,\ldots,p\}$,  $[u_l,v_l]\cap\mbox{supp}(\nu) \neq \emptyset$. If $a=$$u_l$ or $v_l$ are not in ${\rm supp}(\nu)$
that is if ${\rm supp}(\nu)$ does not stick to the frontier of $\mathbb{R} \setminus {\cal O}_{i}$ at these points, then the previous authors established that the density exhibits behavior
closely resembling that of $\sqrt{\vert x-d\vert}$ for $x$ near $d=\phi_i(a)$. We will say that such an edge $\phi_i(a)$ is regular. This is for instance obviously always the case dealing with a discrete measure $\nu$.\\
Nevertheless for some measures $\nu$ with a  density decreasing quite fast to zero at an  edge  of the support of $\nu$, such an edge may coincide with some $u_l$ or $v_l$, that is ${\rm supp}(\nu)$ may stick to the frontier of $\mathbb{R}\setminus  {\cal O}_{i}$ at this point. Then, at the corresponding edge of the support of $\mu_i$, the density $p$ may exhibit different behaviour.  This can be seen for instance in the following example investigated by  Lee and Schnelli \cite{LS} : 
$$d\nu(x) := Z^{-1}(1 + x)^a(1 -x)^b f(x)1_{[-1,1]}(x)\mathrm{d}x$$
where $a<1,b>1$, f is a strictly positive ${\cal C}^1$-function and Z is a normalization constant. 
Indeed let $\sigma_0$ be such that 
$$\int \frac{1}{(1-x)^2}\mathrm{d}\nu(x)=\frac{1}{\sigma_0^2}.$$
Let us consider 
$$\mathbb{R} \setminus {\cal O}_1={\rm supp} (\nu) \cup \{ u \in \mathbb{R} \setminus {\rm supp} (\nu), \int \frac{1}{(u-x)^2}\mathrm{d}\nu(x) \geq \frac{1}{\sigma^2}\}.$$
It can be easily seen that for all
{$\sigma> \sigma_0$}, $\mathbb{R}\setminus {\cal O}_1=[u_\sigma, v_\sigma]$  with $u_\sigma< -1 <1 < v_\sigma,$ so that \begin{equation}\label{regular} \mbox{ ~supp}\left(\mu_{sc}\boxplus \nu\right)=[\phi_1(u_\sigma),\phi_1(v_\sigma)];\end{equation} thus, $\phi_1(v_\sigma)$ is a regular edge and we have ${ p(x)\sim C (\phi_1(v_\sigma)-x )^{\frac{1}{2}}}$.
~~\\ ~~
  Now, for all  {$ \sigma\leq  \sigma_0$}, one can see that 
$\mathbb{R}\setminus  {\cal O}_1=[u_\sigma, 1]$, with
$u_\sigma< -1$, so that $\mbox{supp}\left(\mu_{sc}\boxplus \nu\right)=[\phi_1(u_\sigma),\phi_1(1)]$; it turns out that the density exhibits the following behaviour at the right edge \begin{equation}\label{notregular}{ p(x)\sim C (\phi_1(1)-x )^{b}}.\end{equation}

~~

We illustrate by  the following picture the difference of behaviour of the density $p$ at edges of $\mu_i$ depending on wether the support of $\nu$ sticks to the frontier of $\mathbb{R}\setminus  {\cal O}_i$ or not. We consider a measure $\nu$ whose support has  three connected 
components $[a_i,b_i], i=1,2,3$. Then, we know that $\mathbb{R}\setminus {\cal O}_i$ has at most three connected components and  each of them 
contains at least a connected component of the support of $\nu$. We draw one possible case where $[a_1,b_1]$ and $[a_2,b_2]$ are in the same connected component $[u_1,v_1]$ of $\mathbb{R}\setminus  {\cal O}_i$ and $b_2=v_1$  whereas $[a_3,b_3]$ is in an other  connected component $[u_2,v_2]$ of $\mathbb{R} \setminus {\cal O}_i$. \\
~~~~

 \noindent  {{$~ {\rm supp} (\nu)\subset \mathbb{R}\setminus  {\cal O}_{i} $ }~~ \\
\vspace*{-4.5cm}\\
\hspace*{2cm}
\setlength{\unitlength}{1cm}
\begin{picture}(6,6)(-1,-1)
   \put(-1,0){\vector(1,0){9}}
   \multiput(0,-0.2)(2,0){2}{${[}$} 
   \multiput(0.5,-0.2)(2,0){1}{${]}$} 
\put(3,-0.2){${]}$} 
    {\multiput(-0.5,-0.2)(5,0){1}{${\bigg[}$} 
   \multiput(3,-0.2)(3,0){1}{${\bigg]}$}} 
    \multiput(5.5,-0.2)(2,0){1}{${[}$} 
   \multiput(6.5,-0.2)(2,0){1}{${]}$} 
 {  \multiput(5,-0.2)(5,0){1}{${\bigg[}$} 
   \multiput(7,-0.2)(3,0){1}{${\bigg]}$} }
 {  \put(-0.6,-0.9){${u_1}$}
   \put(3.1,-0.9){${=v_1}$} 
   \put(4.9,-0.9){${u_2}$}
   \put(6.9,-0.9){${v_2}$}}

 \put(0,-0.9){${a_1}$}
   \put(0.5,-0.9){${b_1}$} 
   \put(1.9,-0.9){${a_2}$}
   \put(2.7,-0.9){${b_2}$}
\put(5.5,-0.9){${a_3}$}
   \put(6.4,-0.9){${b_3}$}

 \end{picture}\\

 {{${\rm supp}  (\mu_{i})$}\\
 
\vspace*{-4.5 cm}
\hspace*{1cm}
\setlength{\unitlength}{1cm}
\begin{picture}(6,6)(-1,-1)
   \put(-1,0){\vector(1,0){9}}
   \multiput(0,-0.2)(3,0){1}{${[}$} 
   \multiput(2,-0.2)(3,0){1}{${]}$} 
   \multiput(5,-0.2)(5,0){1}{${[}$} 
   \multiput(7,-0.2)(3,0){1}{${]}$} 
   \put(-0.6,-0.9){${\phi_{i}(u_1)}$}
   \put(1.9,-0.9){${\phi_{i}(v_1)}$} 
   \put(4.9,-0.9){${\phi_{i}(u_2)}$}
   \put(6.9,-0.9){${\phi_{i}(v_2)}$}
 \end{picture}}
~~\\

\hspace*{1,1cm} \hspace*{1cm}$\uparrow$  \hspace*{1.8cm}$\uparrow$ \hspace*{2.9cm}$\nwarrow~~~~~~~ \nearrow$ \\   
\hspace*{2cm}
  \hspace*{-1cm} $p(x)\sim C \vert d-x \vert ^{\frac{1}{2}}$ $\begin{array}{cc} \mbox{the singularity of~} p \\ \mbox{may   change!}
\end{array}$ \hspace*{0.5cm} $p(x)\sim C \vert d-x \vert ^{\frac{1}{2}}$


\section{Outliers of general spiked models}
In Section \ref{seminal}, we presented the seminal works on the behavior of the largest eigenvalues for finite rank deformations of standard models. 
It turns out that the previous analysis in Section 3 allows to understand  the appearence of outliers
of general spikes models that is when  $A_N$ is a deformation with full rank and provides the good way to generalize the pioneering works.
Actually, the relevant criterion for a spiked eigenvalue of $A_N$ to generate an outlier in the spectrum of the deformed model is to belong to some set related to the subordination functions. In the finite rank case this criterion reduces to a critical threshold.\\
 Note that in the iid case,  we do not assume anymore that 
 the entries are   Gaussian; the  results stated in this section  are obtained under different technical  assumptions on the entries on the non-deformed model that we do not precise here and we refer the reader to the corresponding papers.\\
In order to adopt universal notations for the three types of deformations, we set 
$$\tilde{A}_N=\left\{\begin{array}{ll} A_N \mbox{~~~~~~~~for additive or multiplicative deformations~} \\
A_N A_N^* \mbox{~~~~for Information-plus-noise type deformation. ~} \end{array}
\right.$$
Thus, for each type of deformation, we assume the following on the perturbation $\tilde{A}_N$:
\begin{itemize}
\item $\mu_{\tilde A_N}$ weakly converges towards a probability measure $\nu$ whose support is compact.
\item The eigenvalues of $\tilde A_N$ are of two types :
\begin{itemize}
\item $N-r$ (r fixed) eigenvalues $\alpha_i(N)$ such that
$$\max_{i=1}^{N-r} {\rm dist}( \alpha_i(N),{\rm supp}(\nu)) \underset{N \to \infty}{\longrightarrow}  0$$
\item a finite  number $J$ of fixed (independent of $N$) eigenvalues called spikes  $\theta_1 > \ldots> \theta_J$ ($>0$ for multiplicative deformations and information-plus-noise type models), $\forall i=1,\ldots,J$, $\theta_i \not\in {\rm supp}(\nu)$,
each $\theta_j$ having a fixed  multiplicity $k_j$, $\sum_{j}k_j=r$. 
\end{itemize}
\end{itemize} For technical reasons, for information-plus-noise type model, we assume moreover that $A_N$ is a rectangular matrix of the type \begin{equation}\label{diagonale} A_N=\begin{pmatrix}  a_1(N) ~~~~~~~~~~~~~~~~(0)\\
~~~(0)\\ ~~~~~~~~~~\ddots~~~~~~~~~~( 0)\\ ~(0)~~~~~~~~~~~~~~~~~~~~\\  ~~~~~~~~~~~~~a_{N}(N)~~~ ( 0  )    \end{pmatrix} \end{equation} 
\subsection{Location of the outliers} \label{Sec-outliers}

{Here is a naive intuition for general additive deformed models in order to make the reader understand the occurence and role of free subordination functions.}
We have the following subordination property
$$g_{\mu \boxplus \nu }(z)= g_{\nu}(\omega_{\mu,\nu}(z))$$ 
For an Hermitian deformed model such that $$M_N= W_N+A_N; ~~\mu_{W_N}\rightarrow \mu; \mu_{A_N} \rightarrow \nu, \mu_{M_N} \rightarrow \mu\boxplus\nu,$$
 the intuition is that  (see \eqref{deteq} below)
 $$g_{\mu_{M_N}}(z) \approx g_{\mu_{A_N}}(\omega_{\mu,\nu}(z)).$$ 
Assume that $A_N$ has a spiked eigenvalue $\theta$ outside its limiting support.
If  $\rho \notin  {\rm supp} \left(\mu \boxplus \nu\right)$ is a solution of  $\omega_{\mu,\nu}(\rho)=\theta$, 
 $g_{\mu_{M_N}}(\rho) \approx g_{\mu_{A_N}}(\omega_{\mu,\nu}(\rho))$
 explodes! \\
~~\\
Therefore the 
conjecture is 
that the spikes $\theta$'s of the perturbation $A_N$  that may generate outliers in the spectrum of $M_N$  belong to   $\omega_{\mu,\nu}\left(\mathbb{R}\setminus  {\rm supp} \left( \mu \boxplus \nu\right)\right) $ and more precisely that  for large $N$, the $\theta$'s such that  
the equation $$\omega_{\mu,\nu}(\rho)=\theta_i$$ has solutions $\rho$ outside  ${\rm supp} \left( \mu \boxplus \nu\right)$ generate  eigenvalues of $M_N$  in  a neighborhood of each of these $\rho$.

This intuition in fact corresponds to true results for both models: i.i.d and isotropic. Nevertheless, their proofs are different. In the following, we present the distinct approaches.
\subsubsection{The i.i.d case}\label{iidspikes}
In this section, we will denote by $\phi, \varphi, {\cal O}$ any of $\phi_i, \varphi_i, {\cal O}_i$ for $i=1,2,3$ introduced in Section 3.2 related  to the investigated deformed model. We choose to present a unified result covering the three types of deformations. Nevertheless, results concerning the Information-Plus-Noise type model involve some technical  additionnal asumptions. We brievely precise them in remarks  following the unified result and refer the reader to the corresponding papers. \\
~~

\noindent {\bf a) A deterministic equivalent}\label{sectioneqdet} \\
In the three deformed models, a deterministic measure  plays  a central role in the study of the spectrum of the deformed models. This measure is a very good approximation of the spectral measure $\mu_{M_N}$ in the sense that almost surely, for large $N$, each interval in   the complement of the  support of
this deterministic measure contains no eigenvalue of $M_N$. This was first established by Bai and Silverstein \cite{BS98} in the multiplicative case.
We now express  this deterministic measure $\nu_N$ in the three models :
\begin{itemize}
\item[i)] Additive deformation of a Wigner matrix \cite{CDFF}:
$$ \nu_N = \mu_{sc} \boxplus \mu_{A_N} .$$
\item[ii)] Multiplicative deformation of a sample covariance matrix (\cite{BS98})
$$\nu_N = \mu_{\operatorname{MP}}\boxtimes{\mu_{A_N}}.$$
\item[iii)] Information plus noise model (\cite{VLM} in the Gaussian case, \cite{BS12})
$$\nu_N = (\sqrt{\mu_{\operatorname{MP}}}\boxtimes_c \sqrt{\mu_{A_N A_N^*}})^2.$$
\end{itemize}
We denote by $K_N$ the support of $\nu_N$ and for $\epsilon >0$, $(K_N)_\epsilon$ denotes an $\epsilon$ neighborhood of $K_N$. We have the following result :
\begin{proposition}
 $\forall \epsilon >0$, $\forall [a,b] \subset (K_N)_\epsilon$,
 \begin{equation} \label{noeigenvalue}
 \mathbb P (\mbox{ for large $N$, $M_N$ has no eigenvalue in~} [a,b] ) = 1
 \end{equation}
 \end{proposition}
\begin{remark}
For non-Gaussian Information plus noise model,  the result is proved only for $a>0$.
\end{remark}
 \noindent
The proof of \eqref{noeigenvalue} in \cite[Theorem 5.1]{CDFF}  relies on the estimation 
\begin{equation} \label{deteq}
\mathbb{E}\left( g_{\mu_{M_N}} (z) \right) - g_{\nu_N} (z) = \frac{1}{N} L_N(z) + O(\frac{1}{N^2})
\end{equation}
with an explicit formula for $L_N$.
\eqref{deteq} is established using an integration by parts formula in the Gaussian case (and an approximate integration by parts formula in the general Wigner case). Note that  the same method was used in \cite{VLM} for Gaussian Information-Plus-Noise type matrices. \\
To prove  the inclusion of the spectrum \eqref{noeigenvalue} for sample covariance matrices and information-plus-noise type matrices \cite{BS98,BS12}, Bai and Silverstein provide  a different approach although it also makes use of Stieltjes transform. \\

{\bf b) An exact separation phenomenon} \label{sectionseparation} \\
A next step in the analysis of the spectrum of deformed iid models is 
an exact separation phenomenon between the spectrum of $M_N$ and the
spectrum of $\tilde A_N$, involving the subordination functions : 
to a gap in the spectrum of $\tilde A_N$, it corresponds, through the  function $\phi$ defined in Section 3.2,  a gap in the
spectrum of $M_N$ which splits the spectrum of $M_N$ exactly as that of $\tilde A_N$. Let $[a,b]$ be a compact set such that for some $\delta>0$, for all large $N$, 
$[a-\delta,b+\delta] \subset \mathbb{R} \setminus {\rm supp}~ (\nu_N)$.
Then, almost surely, for large $N$, $[\varphi(a), \varphi(b)]$ is in the complement of the spectrum of $\tilde A_N$.
Hence, with the convention that for any $N\times N$ matrix $Z$, $\lambda _0(Z)=+\infty $ and 
$\lambda _{N+1}(Z)=-\infty $, 
there is $i_N\in \{ 0, \ldots , N\}$ such that 
\begin{equation}{\label{sep1}}
\lambda _{i_N+1}(A_N) < \varphi(a)\quad 
\text{and} \quad \lambda _{i_N}(A_N) > \varphi(b). 
\end{equation}

 Moreover,  $[a,b]$ splits the spectrum of $M_N$ exactly as $[\varphi(a),\varphi(b)]$
splits the spectrum of $\tilde A_N$ as stated by the following. 

\begin{theorem}
With $i_N$ satisfying \eqref{sep1}, one has
\begin{equation}\label{sep2}
\mathbb P[\lambda _{i_N+1}(M_N) < a \, \text{ and } \, \lambda _{i_N}(M_N) > b, \, \text{for all large $N$}]=1.\\
\end{equation}
\end{theorem}
\begin{remark}
For non-Gaussian Information plus noise model, if $c<1$, the result is proved only for $b$ such that $\varphi(b)>0$.
\end{remark}
The following picture illustrates this exact separation phenomenon. \\

\noindent 
  $$[a,b] \subset ~\mathbb{R} \setminus   {\rm supp} (\nu_N) \longleftrightarrow [\varphi(a), \varphi(b)]$$

 \hspace*{4.8cm} gap in ~Spect($M_N$) $\longleftrightarrow$
 gap in ~Spect($\tilde A_N$)

\vspace*{-4cm}
\setlength{\unitlength}{1cm}
\hspace*{2cm}\begin{picture}(6,6)(-1,-1)
   \put(-1,0){\vector(1,0){9}}
   \multiput(0,-0.2)(2,0){4}{\line(0,6){0.4}} 
   \put(-1,-0.6){}
   \put(-1.1,-0.6){$\cdots~~\lambda_{l+1}(\tilde A_N)$} 
   \put(1.9,-0.6){$\varphi(a)$}
   \put(4,-0.6){$\varphi(b)$}
   \put(5.9,-0.6){$\lambda_l(\tilde A_N)\cdots$}
 \end{picture}\\ 
\vspace*{-0.7cm}
\hspace*{2cm}  \hspace*{-0.5cm}$\underbrace{\hspace*{3.8cm}}$ \hspace*{2.1cm} $\underbrace{\hspace*{3.5cm}}$\\
\hspace*{2cm} \hspace*{-0.5cm}$\mbox{$N-l$ eigenvalues of} ~\tilde A_N $ \hspace*{2.3cm}$\mbox{$l$ eigenvalues of} ~\tilde A_N $\\
 ~\\

 \vspace*{-4.5cm}
\hspace*{2cm}\begin{picture}(6,6)(-1,-1)
   \put(-1,0){\vector(1,0){9}}
   \multiput(0,-0.2)(2,0){4}{\line(0,6){0.4}} 
   \put(-1,-0.6){}
   \put(-1.1,-0.6){$\cdots ~~ \lambda_{l+1}(M_N)$} 
   \put(1.9,-0.6){$a$}
   \put(4,-0.6){$b$}
   \put(5.9,-0.6){$\lambda_l(M_N)$~~$\cdots$}
 \end{picture}\\ 
\vspace*{-0.7cm}
\hspace*{2cm}  \hspace*{-0.5cm}$\underbrace{\hspace*{3.8cm}}$ \hspace*{2.1cm} $\underbrace{\hspace*{3.5cm}}$\\
\hspace*{2cm} \hspace*{-0.5cm}$\mbox{$N-l$ eigenvalues of} ~M_N$ \hspace*{2.2cm}$\mbox{$l$ eigenvalues of} ~M_N$\\
 ~\\

Again, this was first observed by Bai and Silverstein \cite{BS99}, in the case of sample covariance matrices.
We refer to \cite{CDFF} for deformed Wigner matrices and to Loubaton-Vallet \cite{LV} (Gaussian case), Capitaine \cite{MC2014} for information plus noise type models.
This exact separation phenomenon leads asymptotically to  the relation \eqref{palier} between the cumulative distribution function  of the $\mu_{i}$'s and the cumulative distribution function of $\nu$.\\

{\bf c) Convergence of eigenvalues} \\
The following result gives the precise statement of the intuition given at the beginning of Section \ref{Sec-outliers} and is  a consequence of the inclusion of the spectrum and the exact separation.
\begin{theorem}\label{outliers}\cite{BY, RaoSil09, CDFF, MC2014} Assume that the LSD $\nu$ of $A_N$ has a finite number of connected components.
For each spiked eigenvalue $\theta _j$, 
we denote by \\
$n_{j-1}+1, \ldots , n_{j-1}+k_j$ the descending ranks of $\theta _j$ among the eigenvalues of $\tilde A_N$. \\
With the notations of Section \ref{sectionsupport},
\begin{itemize}
\item[{1)}] If $\theta_j \in {\cal O}$, 
the $k_j$ eigenvalues $(\lambda_{n_{j-1}+i}(M_N), \, 1 \leq i \leq k_j)$ 
converge almost surely outside the support of $\mu $ 
towards $\rho _{\theta _j}=\phi(\theta _j)$.
\item[\text{2)}] If $\theta_j \in \mathbb R \backslash {\cal O}$ 
then we let $[s_{l_j}, t_{l_j}]$ (with $1\leq l_j\leq m$) be the connected component 
of $\mathbb R \backslash {\cal O}$ which contains $\theta _j$.
\begin{itemize}
\item[{a)}] If $\theta_j$ is on the right (resp. on the left) of any connected component of ${\rm supp}(\nu )$ 
which is included in $[s_{l_j},t_{l_j}]$ then the $k_j$ eigenvalues $(\lambda_{n_{j-1}+i}(M_N)$, $1\leq i\leq k_j)$ 
converge almost surely to $\phi(t_{l_j}^+)$ (resp. $\phi(s_{l_j}^-)$) 
which is a boundary point of the support of $\mu $. 
\item[{b)}] If $\theta_j$ is between two connected components of ${\rm supp}(\nu )$ 
which are included in $[s_{l_j},t_{l_j}]$ then 
the $k_j$ eigenvalues $(\lambda_{n_{j-1}+i}(M_N)$, $1\leq i\leq k_j)$ 
converge almost surely to the $\alpha _j$-th quantile of $\mu $ 
(that is to $q_{\alpha _j}$ defined by $\alpha _j=\mu  (]-\infty , q_{\alpha _j}])$) 
where $\alpha _j$ is such that $\alpha _j=1-\lim _N\frac{n_{j-1}}{N}=\nu (]-\infty ,\theta _j])$.
\end{itemize}
\end{itemize}
\end{theorem}
\begin{remark}
For Information-Plus-Noise type models, in 2)a) of Theorem \ref{outliers}, if $\theta_j$ is on the left of any connected component of  ${\rm supp}(\nu )$ 
which are included in  the {\it first} connected component  $[s_1,t_1]$ of $\mathbb R \backslash {\cal O}$, then the convergence of the corresponding eigenvalues of $M_N$ 
towards $\phi(s_{1}^-)$ is established only for $\phi(s_{1}^-)=0$.
\end{remark}

\subsubsection{The isotropic case}
Here we  consider an  additive spiked deformation of an isotropic matrix
$$M_N=A_N+U_N^*B_NU_N.$$ $U_N$ is a unitary matrix whose distribution is  the normalized Haar measure on the unitary group ${\rm U}(N)$.
$B_{N}$ is a  deterministic Hermitian matrix of size
$N\times N$ such that 
 $\mu_{B_{N}}$ converges weakly to $\mu$ compactly supported as $N\to\infty$ and 
such that  the eigenvalues of $B_{N}$
converge uniformly to $\text{supp}(\mu)$ as $N\to\infty$. $A_N$ is a deterministic Hermitian 
$N\times N$ perturbation as defined at the beginning of this Section.\\

Note that if  $A_N$ has no outlier,
that is if $\{\theta_{1},\cdots,\theta_{J}\}=\emptyset$, then the general study of Collins and Male allows to deduce that neither does $M_N$ (see 
Corollary 3.1 in \cite{ColMal11}) meaning that for all large $N$, all  the eigenvalues of $X_N$ are inside a small neighborhood of the support of 
$\mu\boxplus \nu$. \\

Assume now that $A_N$ has outliers. 
Then Belinschi, Bercovici, Capitaine and F\'evrier established in \cite{BBCF} the following result.
\begin{theorem}\label{Main+}
Set $K=\mathrm{supp}(\mu\boxplus\nu)$,
$$K'=K\cup \left[\bigcup_{i=1}^J
\omega_2^{-1}(\{\theta_i\})\right],$$
whith $\omega_2$  defined as in Section 2.2.1.
The following results hold almost surely for large $N${\rm:}

\noindent  Given $\varepsilon>0$ and denoting by $K'_\epsilon$ an $\epsilon$ neighborhood of $K$, we have $${\rm spect}(M_N)\subset K'_\varepsilon.$$
\noindent Let $\rho$ be a  fixed number in $ K'\setminus K$  and $\theta_i $ be such that $\omega_2(\rho)=\theta_i$. For any $\varepsilon>0$  such that $(\rho-2\varepsilon,
\rho+2\varepsilon)\cap K'=\{\rho\}$, we have
$$
\mathrm{card}(\{{\rm spect}(M_N)\cap(\rho-\varepsilon,\rho+\varepsilon)\})=k_i.
$$

\end{theorem}

Here we explain the sketch of the proof. Fix $\alpha\in\text{supp}(\nu)$.
Due to the left and right invariance of the Haar measure on ${\rm U}(N)$ we may assume without loss of 
generality that both $A_{N}$ and $B_{N}$ are diagonal matrices. More precisely, we let $A_N$ be the diagonal matrix
$$
A_N=\text{ Diag}
(\underbrace{\theta_1,\dots,\theta_1}_{k_1 \mbox{times}},\dots,\underbrace{\theta_J,\dots,\theta_J}_{k_J \mbox{times}},\alpha_1^{(N)},\dots,\alpha_{N-r}^{(N)}),
$$
and  write  $A_{N}=A_{N}'+A_{N}'',$ 
where 
\[
A_{N}'=\text{Diag}(\underbrace{\alpha,\ldots,\alpha}_r,\alpha_{1}^{(N)},\ldots,\alpha_{N-r}^{(N)}),
\]
 and
\[
A_{N}''= {\text{Diag}(\underbrace{\theta_1-\alpha,\dots,\theta_1-\alpha}_{k_1 \mbox{times}},\dots,\underbrace{\theta_J-\alpha,\dots,\theta_J-\alpha}_{k_J \mbox{times}},\underbrace{0,\ldots,0}_{N-r})}.
\]
We have
$A_{N}''=P_N^*\Theta P_N$, where $P_N$ is the $r\times N$ matrix representing the usual projection 
$\mathbb{C}^N\to\mathbb{C}^r$ onto the first $r$ coordinates, and 
\[
\Theta= {\text{Diag}(\underbrace{\theta_1-\alpha,\dots,\theta_1-\alpha}_{k_1 \mbox{times}},\dots,\underbrace{\theta_J-\alpha,\dots,\theta_J-\alpha}_{k_J \mbox{times}}}).
\]

 The matrices $A_{N}'$ and $B_{N}$ have no spikes, 
and therefore \cite[Corollary 3.1]{ColMal11}
applies to the matrix $M'_N=A_{N}'+U_{N}^{*}B_{N}U_{N}$. Note that the limiting spectral measure is still $\mu\boxplus \nu$.
The  first key idea is due to Benaych-Georges and Nadakuditi  \cite{BN1} and consists in reducing  
the problem of locating outliers of the deformations to a  convergence problem of a 
 fixed size $r\times r$ random matrix, by  using the Sylvester's determinant identity: for rectangular matrices $X$ and $Y$ such that  $XY$
and $YX$ are square, we have
\begin{equation}\label{Syl}
\det(I+XY)=\det(I+YX).
\end{equation}
Given $z$ outside the support of $\mu\boxplus \nu$, we have
\[
\det(zI_{N}-(A_{N}+U_{N}^{*}B_{N}U_{N}))=\det(zI_{N}-M_{N}')\det(I_{N}-\left(zI_{N}-M_{N}'\right)^{-1}\,P_N^*\Theta P_N).
\]
 so that using Sylvester's identity, 
\[
\det(zI_{N}-(A_{N}+U_{N}^{*}B_{N}U_{N}))=\det(zI_{N}-M_{N}')\det(I_{r}-P_N\left(zI_{N}-M_{N}'\right)^{-1}\,P_N^*\Theta).
\]
 We conclude that  the eigenvalues of $A_{N}+U_{N}^{*}B_{N}U_{N}$
outside $\mu\boxplus \nu$ are precisely the zeros of the function $\det(F_N(z))$, where
\begin{equation}
F_{N}(z)=I_{r}-P_N\left(zI_{N}-M_{N}'\right)^{-1}\,P_N^*\Theta\label{MN}.
\end{equation}

The key idea  is now  to establish an approximate matricial subordination result. Biane \cite{Biane98} proved the 
stronger result 
that for any  ${\bf a}$ and ${\bf b}$  free selfadjoint random variables in a tracial W*-probability space,
 there exists an analytic self-map $\omega\colon\mathbb
C^+\to\mathbb C^+$ of the upper half-plane   so that 
\begin{equation}\label{BianeCondExp}
\mathbb E_{\mathbb C[{\bf a}]}\left[(z-({\bf a+b}))^{-1}\right]=(\omega(z)-{\bf a})^{-1},
\quad z\in\mathbb C^+.
\end{equation}
Here $\mathbb E_{\mathbb C[{\bf a}]}$ denotes the conditional expectation onto the
von Neumann algebra generated by ${\bf a}$.  It can be proved that an approximate version does hold in the sense that 
the compression 
$$ P_N \left[\mathbb E\left[(z-(A_N'+U_N^*B_NU_N))^{-1}\right]^{-1}+A_N' \right] P_N
$$
 is close to $\omega_2(z)I_{r}$,  as $N$ goes to infinity.
Thus, it turns out that almost surely  the sequence  $\{F_N\}_N$ converges uniformly on compact subsets of $\mathbb{C}\setminus \rm{supp}(\mu\boxplus\nu)$ to the analytic function $F$ defined by
\[
F(z)=\mathrm{Diag}\left(\underbrace{1-\frac{\theta_{1}-\alpha}{\omega_2(z)-\alpha}}_{k_1 \mbox{times}},
\dots, \underbrace{1-\frac{\theta_{J}-\alpha}{\omega_2(z)-\alpha}}_{k_J \mbox{times}}\right),
\]
where $\omega_2$ is the subordination function from (\ref{subord1}).
The set of points $z$
such that $F(z)$ is not invertible is precisely $\bigcup_{i=1}^J\omega_2^{-1}(\{
\theta_i\})$.\\

It follows from this result that a remarkable new phenomenon arises: a 
single spike of $A_{N}$ can generate asymptotically a finite or even a countably infinite
set of outliers of $M_{N}$. This arises from the fact that the restriction
to the real line of some subordination functions may be many-to-one, that is, with the  above
notation, the set $\omega_i^{-1}(\{\theta\})$ may have cardinality strictly
greater than 1, unlike the subordination function related to free convolution with
a semicircular distribution that was used in section \ref{iidspikes}.
The following numerical simulation, due to Charles Bordenave, illustrates the appearance of two outliers 
arising from a single spike.  We take $N=1000$ and $M_N=A_N+
U_NB_NU_N^*$, where $B_N={\rm{ diag
}}(\underbrace{-1,\ldots
  ,-1}_{\frac{N}{2}}, \underbrace{1,\ldots,
 1 }_{\frac{N}{2}})$, and
$$
A_N=\left[\begin{array}{cc}
\frac{W_{N-1}^G}{2} & 0_{(N-1)\times 1}\\
0_{1\times (N-1)} & 10
\end{array}\right],
$$
with $W_{N-1}^G$ being sampled from a standard $999\times999$ G.U.E..
{This is not a spiked deformed GUE model and now,  the spike $\theta=10$ is associated to
the matrix approximating the semicircular distribution.}
We have the subordination identities 
$$g_{\mu_{sc} \boxplus \nu }(z)= g_{\nu}(\omega_2(z))= g_{\mu_{sc}}(\omega_1(z))$$
where $\omega_2$ is injective on $\mathbb{R}\setminus  {\rm supp} \left(\mu_{sc} \boxplus \nu\right)$ but $\omega_1$ may be many to one. Actually 
 $\omega_1(\rho)=10$ has  2 solutions $\rho_1$ and $\rho_2$.
\noindent\begin{center}
     \includegraphics[width=12cm,height=7cm]{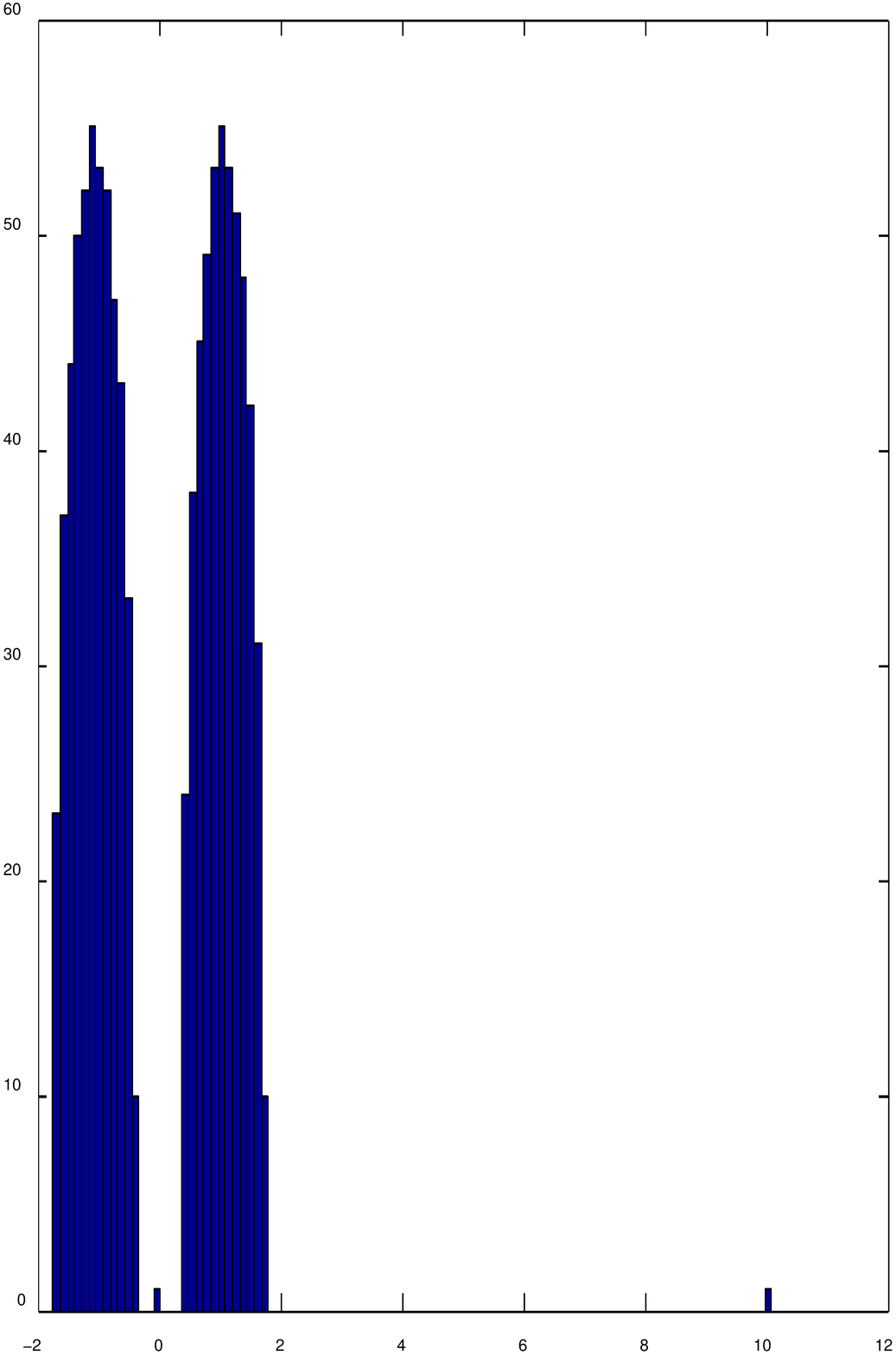}     
\par\end{center}

Actually, one may consider additive models where both $A_N$ and $B_N$ have spiked outliers (\cite{BBCF}).
Let us consider 
$$\hspace*{-0.2cm}M_N =U_NB_N U_N^* + A_N$$ where $U_N $ is a Haar unitary matrix, $ A_N$ and $B_N$ are deterministic diagonal matrices
such that $\mu_{B_N} \underset{N\to +\infty}{\longrightarrow} \mu$
and $\mu_{A_N} \underset{N\to +\infty}{\longrightarrow} \nu$, both $\mu$ and $\nu$ being compactly supported. We also assume that there exist  a spiked
$\theta \notin {\rm supp}(\nu)$ which is an eigenvalue of $A_N$  with multiplicity $k$
and a spiked $\alpha \notin {\rm supp}(\mu)$ which is an eigenvalue of $B_N$ with multiplicity $l$,
whereas the other eigenvalues are uniformly close to the limiting supports.
We have the subordination properties
$$g_{\mu \boxplus \nu}(z)=g_\nu(\omega_{2}(z))=g_\mu(\omega
_{1}(z)).$$
{If there exists $\rho \in \mathbb{R} \setminus {\rm supp}(\mu \boxplus \nu)$ such that $$\left\{\begin{array}{ll}\omega_{1}(\rho)= \alpha\\\omega_{2}(\rho)=\theta \end{array} \right.$$
then for all large $N$, there are $k+l$ outliers of $M_N$ in a neighborhood of $\rho$.}\\

Such results are established for
 multiplicative perturbations of unitarily invariant matricial models, based on similar ideas, with
the subordination function  replaced by its multiplicative counterpart.
\subsection{Eigenvectors}
 For a general perturbation, dealing with sample covariance matrices, S. P\'ech\'e and O. Ledoit \cite{PL} introduced
a tool to study the average behaviour of the eigenvectors but it seems that this
did not allow them to focus on the eigenvectors associated with the eigenvalues
that separate from the bulk.
It turns out that further studies \cite{MCJTP, BBCF, C2015} point out that
the angle between the eigenvectors of the outliers of the deformed model  and the eigenvectors associated to the corresponding  original spikes
is determined by free subordination functions.

The following theorem \cite{MCJTP, BBCF} holds for spiked additive deformations  in the i.i.d case as well as in the isotropic case. 
Let $A_N$ be a deterministic  deformation as defined as the beginning of Section  4; dealing with either a deformed Wigner matrix or a deformed isotropic matrix as defined at the beginning of Section 4.1.2, we have the following
\begin{theorem}
Set $K=\mathrm{supp}(\mu\boxplus\nu)$,
$$K'=K\cup \left[\bigcup_{i=1}^p
\omega_2^{-1}(\{\theta_i\})\right],
$$
and let $\omega_2$ be the subordination function satisfying {\rm (\ref{subord1})}.
 Let $\rho$ be a  fixed number in $ K'\setminus K$  and $\theta_i $ be such that $\omega_2(\rho)=\theta_j$.  Let $\varepsilon>0$ be such that $(\rho-2\varepsilon,
\rho+2\varepsilon)\cap K'=\{\rho\}$. Let $\xi$ be a unit eigenvector associated to an eigenvalue of $M_N$ in $(\rho-\epsilon,\rho+\epsilon)$. Then when $N$ goes to infinity,
$$
\Vert P_{\operatorname{Ker~}(\theta_l I_N -A_N)} (\xi)\Vert^2 \rightarrow \frac{\delta_{jl}}{\omega'_2(\rho)}\mbox{~~almost surely }$$

\end{theorem}

Similar results are established for spiked multiplicative deformations  in the i.i.d case as well as in the isotropic case \cite{MCJTP, BBCF}
and for information-plus-noise type models in \cite{C2015} in the i.i.d case. 
See the following tabular.
Note that in the i.i.d case everything is explicit and can  be rewritten as follows using  the $\phi_i$'s defined in Section 3.2.
\begin{theorem}
Let $n_{j-1}+1, \ldots , n_{j-1}+k_j$  the descending ranks of $\theta _j$ among the eigenvalues of $\tilde A_N$ and  
 $\xi(j)$ a unit eigenvector associated to  one of the eigenvalues $(\lambda_{n_{j-1}+q}(M_N)$, $1\leq q\leq k_j)$.\\
Then 
when $N$ goes to infinity,
\begin{itemize}
\item  For any $\theta_l\neq \theta_j$,
$${\Vert 
P_{\operatorname{Ker~}(\theta_l I_N -\tilde A_N)} (\xi(j))\Vert \rightarrow 0} \mbox{~~almost surely }$$
\item $${\Vert 
P_{\operatorname{Ker~}(\theta_j I_N -\tilde A_N)} (\xi(j))\Vert^2 \rightarrow \alpha_j }\mbox{~~almost surely }$$

{where $\alpha_j=\left\{\begin{array}{lllll}\phi'_1(\theta_j)=1 -\sigma^2 \int \frac{1}{(\theta_j-x)^2}\mathrm{d}\nu(x) { \mbox{~for  deformed Wigner matrices~}} \\ \\
\frac{\theta_j \phi'_2(\theta_j)}{\phi_2(\theta_j)}=
\frac{1- c\int \frac{x^2}{(\theta_j-x)^2}\mathrm{d}\nu(x)}{1+c \int \frac{x}{(\theta_j-x)}\mathrm{d}\nu(x)} {\mbox{~for sample covariance matrices}}\\ \\
\frac{ \phi_{3}'({\theta_j})}{1+ \sigma^2 cg_\nu(\theta_j)} { \mbox{~for  information-plus-noise type matrices~}} 
 \end{array}\right.
$}

\end{itemize}
\end{theorem}

Here are the common basic ideas of the proof of these results
\cite{MCJTP}.\\
First note that if $u_1,\ldots, u_N$ and $w_1,\ldots, w_N$ are respectively a basis of eigenvectors associated with $\lambda_1(\tilde A_N), \ldots, \lambda_N(\tilde A_N)$
  and  with $\lambda_1(M_N), \ldots, \lambda_N(M_N)$, we have
$$\Tr \left[h(M_N) f(\tilde A_N)\right] =\sum_{k,l} h(\lambda_k(M_N)) f(\lambda_l(\tilde A_N)) \vert \langle u_l,w_k \rangle \vert^2.$$
Thus, since the $\theta_l$'s  separate from the rest of the spectrum of $\tilde A_N$ and the  outliers of $M_N$ separate from the rest of the spectrum of $M_N$, one can deduce the asymptotic norm of the projection onto an eigenspace associated to a spike $\theta_i$, of an eigenvector associated to an outlier of $M_N$ from the 
    study of the asymptotic behaviour of $\Tr \left[h(M_N) f(\tilde A_N)\right] $ for a fit choice of $h$ and $f$.
Then, 
a concentration of measure phenomenon  reduces the problem to the study of $\mathbb{E}( \Tr \left[h(M_N) f(\tilde A_N)\right]) $.\\
The third key point is to approximate   the function $h$ by its convolution by the Poisson Kernel in order to exhibit the resolvent of the deformed model
$$\mathbb{E}\left[ \Tr \left[ h (M_N) f (\tilde A_N)\right] \right]=- \lim_{y\rightarrow 0^{+}}\frac{1}{\pi} \Im \int  \mathbb{E} \left( \Tr \left[G_{N}(t+iy)f (\tilde A_N)\right] \right)h(t) \mathrm{d}t $$
\noindent where $G_N(z) =(zI_N -M_N)^{-1}$.
Finally, writing 
$\tilde A_N= U^*DU$, with $D$ diagonal and $U$ unitary, defining 
 $\tilde{G}_N := UG_NU^*$, the result follows from sharp estimations of~$\mathbb{E} (\{\tilde{ G}_N\}_{kk}(z)) $,
for any $k$ in $\{1, \ldots,N\}$.
\subsection{Unified understanding}
 In conclusion, solving the problem of outliers consists in solving an equation involving the free subordination function and the spikes of the perturbation. Moreover, the norm of the orthogonal projection of an eigenvector associated to an outlier of the deformed model onto the eigenspace of the corresponding spike of the perturbation is asymptotically determined by the  free subordination function. This is summarized in the following tabular. \\
In the tabular, $Y_N$ denotes  a  Hermitian random matrix of iid type ($Y_N = W_N$, $S_N$ or $\sigma \frac{X_N}{\sqrt{p}}$ according to the deformations, see section 1.3)  or $Y_N$ is  unitarily invariant (resp. biunitarily invariant for the information plus noise model).

\vspace{.3cm}
\noindent \begin{tabular}{|l|l|l|} \hline \hspace*{0.3cm}
 $\begin{array}{llllll} ~~\\M_N=A_N+Y_N \\ \mu_{A_N}\rightarrow_{N \rightarrow +\infty}\nu \\\mu_{Y_N}\rightarrow_{N \rightarrow +\infty} \mu \\ {\theta \in \mbox{ Spect}(A_N)}\\{\theta \mbox{ multiplicity~} k_i} \\ { \theta \notin {\rm supp}(\nu)}  \end{array}$  \hspace*{-0.8cm} &$\begin{array}{llllll} ~~\\  \hspace{-0.2cm} M_N=A_N^{1/2}Y_NA_N^{1/2}\\\mu_{A_NA_N^*}\rightarrow_{N \rightarrow +\infty}\nu \\\mu_{Y_N}\rightarrow_{N \rightarrow +\infty}\mu\\{ \theta \in \mbox{ Spect}(A_N)}\\{\theta \mbox{ multiplicity~} k_i }\\ {\theta > 0, \theta \notin {\rm supp}(\nu) }\end{array}$  \hspace*{-0.4cm}&$\begin{array}{lllllll}  ~~\\\hspace*{-0.4cm}M_N=  (A_N+Y_N)(A_N+Y_N)^* \\\mu_{A_NA_N^*}\rightarrow_{N \rightarrow +\infty}\nu \\  \mu_{Y_NY_N^*} \rightarrow_{N \rightarrow +\infty}\mu\\ \sqrt{\mu} \rm{~or~} \sqrt{\nu}~ \boxplus_c \mbox{\small infinitely  divisible}  \\{ \theta \in \mbox{ Spect}(A_NA_N^*)}\\ {\theta \mbox{ multiplicity~} k_i} \\ {\theta >0, \theta \notin {\rm supp}(\nu)} \end{array}$ \\
 \hline \hspace*{0.2cm}$
\mu_{M_N}\rightarrow_{N \rightarrow +\infty} \mu \boxplus{\nu}$  \hspace*{-0.5cm}&$ \hspace{-0.2cm}\mu_{M_N}\rightarrow_{N \rightarrow +\infty} \mu \boxtimes{\nu}  $  \hspace*{-0.5cm}& $
\hspace{-0.2cm} \mu_{M_N}\rightarrow_{N \rightarrow +\infty} (\sqrt{\mu} \boxplus_c \sqrt{\nu})^2 $
 \\ \hline  $\begin{array}{lll} ~~\\ g_\tau(z) = \int_\mathbb{R} \frac{d\tau(x)}{z-x} \\~~ \end{array}$ \hspace*{-0.5cm}& $\begin{array}{lll} ~~\\  \hspace{-0.2cm}\Psi_\tau (z)=\frac{1}{z}g_\tau(\frac{1}{z}) -1\\~~ \end{array}$  \hspace*{-0.5cm} & $\begin{array}{lll} ~~\\  \hspace{-0.4cm} H_{\hspace{-0.05cm}\sqrt{\tau}}^{(c)} \hspace{-0.15cm}=\frac{c}{z} g_\tau(\frac{1}{z})^2 +(1-c) g_\tau(\frac{1}{z}) \\~~ \end{array}$
\\ \hline $\begin{array}{lll} ~~\\g_{\mu \boxplus \nu }(z)= g_{\nu}(\omega_{\mu,\nu}(z))\\~~ \end{array}$ \hspace*{-0.5cm}& $ \Psi_{\mu\boxtimes \nu}(z) =\Psi_\nu(F_{\mu,\nu}(z))$  \hspace*{-0.3cm}& $  \hspace{-0.2cm}
H^{(c)}_{\sqrt{\mu} \boxplus_c \sqrt{\nu}}(z)=H^{(c)}_{\sqrt{\nu}}(\Omega_{\mu,\nu}(z)) $\\ \hline ${ \begin{array}{lllc} k_i \mbox{~outliers  of~} M_N \\\mbox{ \small  in the neighborhood}\\ \mbox{ \small of each $\rho$ s.t}\\
\omega_{\mu,\nu}(\rho)=\theta\end{array}}$  \hspace*{-0.5cm}&${\begin{array}{lllc} k_i \mbox{~outliers  of~} M_N \\\mbox{ \small  in the neighborhood}\\ \mbox{ \small of each $\rho$ s.t}\\ \frac{1}{F_{\mu,\nu}(1/\rho)} 
=\theta\\~~ \end{array}}$  \hspace*{-0.5cm}&${\begin{array}{lllc} k_i \mbox{~outliers  of ~}M_N \\\mbox{ \small  in the neighborhood}\\ \mbox{ \small of each $\rho$ s.t}\\ \frac{1}{\Omega_{\mu,\nu}(1/\rho)}=\theta\\~~ \end{array}}$  \\ \hline
 ${ \begin{array}{lllc} \xi \mbox{~eigenvector of~} M_N \\\mbox{ \small associated to an outlier}\\ \mbox{ \small  in the neighborhood }\\ \mbox{ \small of  $\rho$ s.t $
\omega_{\mu,\nu}(\rho)=\theta$ }\\ \hspace*{-0.4cm} \mbox{ \small $ \Vert P_{\mbox{Ker} (\theta I-A)} \xi \Vert^2 \rightarrow_{N \rightarrow +\infty} \frac{1}{\omega_{\mu,\nu}^{'}(\rho)}$} \end{array}}$  \hspace*{-0.5cm}& ${ \begin{array}{lllc} \xi \mbox{~eigenvector of~} M_N \\\mbox{ \small associated to an outlier}\\ \mbox{ \small  in the neighborhood }\\ \mbox{ \small of  $\rho$ s.t $
\frac{1}{F_{\mu,\nu}(1/\rho)} 
=\theta$ }\\ \hspace*{-0.4cm} \mbox{ \small $ \Vert P_{\mbox{Ker} (\theta I-A)} \xi \Vert^2 \rightarrow_{N \rightarrow +\infty} \frac{\rho F_{\mu,\nu}(1/\rho)}{F_{\mu,\nu}^{'}(1/\rho)}$} \end{array}}$  \hspace*{-0.5cm}&${ \begin{array}{lllc} \xi \mbox{~eigenvector of~} M_N \\\mbox{ \small associated to an outlier}\\ \mbox{ \small  in the neighborhood }\\ \mbox{ \small of  $\rho$ s.t $
\frac{1}{\Omega_{\mu,\nu}(1/\rho)}=\theta$ }\\ \hspace*{-0.4cm} \mbox{ \small $  \Vert P_{\mbox{Ker} (\theta I-A)} \xi \Vert^2 \rightarrow_{N \rightarrow +\infty} \frac{\rho^2 g_{(\sqrt{\mu} \boxplus_c \sqrt{\nu})^2}(\rho)}{\theta^2 g_\nu(\theta) \Omega_{\mu,\nu}^{'}(1/\rho)} $} \end{array}}$   \\ \hline
\end{tabular}
~~~~\\

Note that up to now, the formula in the lower right corner of the previous tabular, concerning  the limiting projection of the eigenvectors associated to outliers of Information-Plus-Noise type models,   has been proved only in the iid case for diagonal perturbation $A_N$  and in the isotropic case for finite rank perturbation $A_N$.
\section{Fluctuations at edges  of spiked deformed models}
In this section, we present results on fluctuations of outliers and eigenvalues at soft edges of the limiting support, with particular stress on understanding the phenomena through free probability theory.

\subsection{The Gaussian case}
\subsubsection{Additive deformation}
Let us consider the deformed G.U.E.. It is known from Johansson \cite{Johansson} (see also \cite{BrezinHikami}) that the joint eigenvalue density induced by the  latter model can be explicitely computed. Furthermore it induces a so-called ``determinantal random point field".\\
When $A_N$ is of finite rank, P\'ech\'e \cite{Peche} obtained a striking phase transition phenomenon for the fluctuations of the largest eigenvalue of the deformed G.U.E..
\begin{theorem} \label{theoTCLpeche}
\begin{enumerate}
\item If $\theta_1 < \sigma$, $\sigma^{-1} N^{2/3}  (\lambda_1(M_N) - 2 \sigma)$ converges in distribution to the G.U.E. Tracy Widom distribution $F_2$.
\item If $\theta_1 = \sigma$, $\sigma^{-1} N^{2/3}  (\lambda_1(M_N) - 2 \sigma)$ converges in distribution to a "generalized" Tracy Widom distribution $F_{3, k_1}$.
\item If $\theta_1 > \sigma$, $ N^{1/2}  (\lambda_1(M_N) - 2 \sigma)$ converges in distribution to the largest eigenvalue of a G.U.E. matrix of size $k_1$ and parameter  $\sigma_{\theta_1} = \sigma \sqrt{1-(\sigma/\theta_1)^2}$. In particular, if $k_1 =1$, $ N^{1/2}  (\lambda_1(M_N) - 2 \sigma)$ converges in distribution to a centered normal distribution with variance $\sigma_{\theta_1}$.
\end{enumerate}
\end{theorem}
The proof is based on the expression of the distribution of the largest eigenvalue in terms of Fredholm determinant and then as a contour integral.
The asymptotic properties  rely on a saddle point analysis. \\

The seminal works concerning full rank deformations of a G.U.E. matrix
made strong assumptions  on the rate of convergence of $\mu_{A_N}$ to $\nu$.
In \cite{Shcherbina2}, the author  investigates  the local edge regime which deals with the behavior of the eigenvalues near any regular extremity point $u_0$ of a connected component of $\text{supp}(\mu_{sc}\boxplus\nu)$. The typical size of the fluctuations of the eigenvalues at  regular edges (see Section  3.2.1) is $N^{-2/3}$. \cite{Shcherbina2} considers the case where $\mu_{A_N}$ concentrate quite fast to the measure $\nu$. In particular, there are no spike. More precisely  
\cite{Shcherbina2} makes a technical assumption on the uniform convergence of the Stieltjes transform of $\mu_{A_N}$ to $g_{\nu}$: 
\begin{equation}\label{Scherbina}\sup_{z \in K} |g_{\mu_{A_N}}(z)-g_{\nu }(z)|\leq N^{-2/3-\epsilon}, \end{equation} where $K$ is some compact subset of the complex plane at a positive distance of the support of $\nu.$ 
Then, 
\cite{Shcherbina2} proves  that the joint distribution of the  eigenvalues converging to $u_0$ have universal asymptotic behavior, characterized by the  Tracy-Widom distribution.
In 
\cite{Adleretal1} and \cite{Adleretal2}, \cite{BK,ABK}
  the authors consider the case where  
$A_N$
has two distinct eigenvalues $\pm a$
of equal multiplicity. They proved the Tracy-Widom  fluctuations at  edges (which are all regular since $\nu$ is discrete).

It turns out that the above  strong assumptions   made on the rate of convergence of $\mu_{A_N}$ to $\nu$ can be removed by  studying the asymptotic distribution of eigenvalues in the vicinity of mobile edges namely the edges of  the deterministic equivalent $\mu_{sc}\boxplus \mu_{A_N}$ of the empirical eigenvalue distribution of the deformed GUE. In \cite{MirSand}, the authors establish the following results. \\

Let $d$ be a regular right edge of $\text{supp}(\mu_{sc}\boxplus \nu)$.
Assume moreover that for any $\theta_j$ such that $\int \frac{\mathrm{d}\nu(s)}{(\theta_j -s)^2} =1/\sigma^2$, we have $d \neq \phi_1(\theta_j)=\theta_j -\sigma^2g_\nu (\theta_j)$.
 It turns out that 
for $\eta$ small enough, for all large $N$, there exists a  unique right edge $d_N$ of  $\text{supp}(\mu_{sc}\boxplus \mu_{A_N})$ in
$]d-\eta, d+\eta[$.
and  the asymptotic distribution of eigenvalues in the vicinity of $d_N$ is universal as the following Theorem \ref{theo: A} states.

\begin{theorem} \label{theo: A}
Let $k$ be a given fixed integer. Let $\lambda_{max}\geq \lambda_{max-1}\geq \cdots \lambda_{max-k+1}$ denote the $k$ largest of those eigenvalues of $M_N$  converging to $d.$
There exists $\alpha >0$ depending on $d_N$ only such that 
 the vector $$\frac{N^{2/3}}{\alpha} \left( \lambda_{max}-d_N, \lambda_{max-1}-d_N , \ldots, \lambda_{max-k+1}-d_N  \right)$$ converges in distribution as $N \to \infty$ to the so-called Tracy-Widom G.U.E. distribution for the $k$ largest eigenvalues (see \cite{TW}).
\end{theorem}

We now turn to the behavior of outliers. Let $\theta_{i}$ be a spiked eigenvalue with multiplicity $k_i$, such that 
$\int \frac{1}{(\theta_i-x)^2}\mathrm{d}\nu(x) <1/\sigma^2$. Recall that in \cite{CDFF}, the authors prove that the spectrum of $M_N$ exhibits $k_i$ eigenvalues in a neighborhood of \begin{equation} \label{defrho} \rho_{\theta_{i}}=\theta_i+\sigma^2\int\frac{\mathrm{d}\nu(x)}{\theta_i-x}.\end{equation}
Oncemore, dealing with mobile edges   related to $\mu_{sc}\boxplus \mu_{A_N}$, \cite{MirSand} obtains the following universal result.

\begin{theorem}
 Let $\theta_i$ be such that $\int \frac{\mathrm{d}\nu(x)}{(\theta_i-x)^2} <1/\sigma^2$ and $\rho_{\theta_i}=\phi_1(\theta_i)$. Then, for $\epsilon>0$ small enough, for all large $N$, {$\text{supp}(\mu_{sc}\boxplus\mu_{A_N})$ has a unique connected component $[L_i(N), D_i(N)]$ inside $]\rho_{\theta_i} -\epsilon,  \rho_{\theta_i} +\epsilon[$}.
Moreover, 
{ the $k_i$ outliers of $M_N$ close to $ \rho_{\theta_i}$  fluctuate at rate $\sqrt{N}$ around $\frac{L_i(N)+D_i(N)}{2} $ as the eigenvalues of a $k_i\times k_i$ GUE.}
\end{theorem}
 The basic tool is a saddle point analysis of the correlation functions of the deformed G.U.E., the contours involving  the image  of $\mathbb{R}$ by the continuous extension of  the subordination  function $\omega_2$ defined by  (\ref{subord1}) with $\mu=\mu_{sc}$ and $\nu=\mu_{A_N}$. 

~

\subsubsection{Sample covariance matrices}
As in the above section, the distribution of the eigenvalues is explicit, with a determinantal structure. The analysis of the fluctuations relies on an expression 
of the distribution of the extremal eigenvalues in terms of a Fredholm determinant and then an asymptotic analysis based on a saddle point method, or a steepest descent method. \\ 
The first result was obtained by Baik,  Ben Arous  and P\'ech\'e \cite{BBP} who described the fluctuations of the largest eigenvalues at the right edge and revealed the  phase transition phenomenon, in the case of a finite rank perturbation $A_N$ of the identity. We refer to \cite{BBP} for the precise statement, their result being an analogue of Theorem \ref{theoTCLpeche}. \\
The full rank case was investigated by Hachem, Hardy and Najim \cite{HHN} for extremal eigenvalues sticking to the bulk (the support of $\mu_{\operatorname{MP}} \boxtimes \nu$). As in Theorem \ref{theo: A}, around regular soft edges,  the  associated
extremal eigenvalues, properly rescaled,  in the vicinity of mobile edges namely the edges of  the deterministic equivalent $\mu_{\operatorname{MP}}\boxtimes \mu_{A_N}$, converge in law to the Tracy-Widom distribution
at the scale $N^{2/3}$.

\subsubsection{Random perturbations}

If one let the perturbation matrix $A_N$  be random then the mobile edges of the equivalent measure become random and may  lead  to different rates of convergence and different asymptotic distributions. We present two examples established respectively by Johansson \cite{J} and Lee and Schnelli \cite{LeeSchnelli} that we revisited through free convolutions.

\begin{itemize}
 \item Johansson \cite{J} considered $$M_N=W_N^G + A_N$$ where $W_N^G$ is a G.U.E. matrix as defined in  Section 1.1.1 and $$A_N= N^{-1/6} \rm{diag}(y_1,\ldots,y_N)$$
where the $y_i$'s are  iid real random variables with distribution $\tau$ and independent from $W_N^G$.
Let us assume that $\tau$ is compactly supported and set $v^2=\int x^2 \mathrm{d}\tau(x)$. Note that almost surely $\mu_{A_N}$ converges weakly to $\delta_0$ and $\mu_{M_N}$ converges weakly to $\mu_{sc}$.
Denote by $\tilde{d}_N$ the {\it deterministic} upper right edge of $\mu_{sc}\boxplus \frac{\tau}{N^{1/6}}$ where $\frac{\tau}{N^{1/6}}$ denotes the pushforward of $\tau$ by the map $x \mapsto \frac{x}{N^{1/6}}$.
 Johansson established that 
\begin{equation}\label{fluctuJohansson} \sigma^{-1}N^{2/3}\left(\lambda_{\rm{max}}(M_N) -\tilde d_N\right) \xrightarrow{\cal D} X+Y \end{equation}
where $X$ and $Y$ are independent random variables, $X$ has the Tracy-Widom distribution and $Y$ has distribution $N(0, \frac{v^2}{\sigma^2})$.
Note that the upper right edge $\tilde d_N$ of $\mu_{sc}\boxplus \frac{\tau}{N^{1/6}}$ is defined by \begin{equation}\label{tilded}\tilde d_N= \tilde t_N + \sigma^2 \int 
\frac{1}{\tilde t_N-x/N^{1/6}}\mathrm{d}\tau(x),\end{equation}
where $\tilde t_N$ in the vicinity of $\sigma$  satisfies \begin{equation}\label{tildet} \int 
\frac{1}{(\tilde t_N-x/N^{1/6})^2}\mathrm{d}\tau(x) =\frac{1}{\sigma^2}.\end{equation}
Consider now the {\it random} upper right edge $ d_N$ of $\mu_{sc}\boxplus \mu_{A_N}$. It is defined by 
\begin{equation}\label{d}d_N=t_N +\sigma^2 \frac{1}{N}\sum_{i=1}^N \frac{1}{t_N -y_i/N^{1/6}},\end{equation}
where $t_N$  in the vicinity of $\sigma$ satisfies 
\begin{equation}\label{t}\frac{1}{N}\sum_{i=1}^N \frac{1}{(t_N -y_i/N^{1/6})^2}=\frac{1}{\sigma^2}.\end{equation}
It is easy to see that (\ref{tilded}), (\ref{tildet}),  (\ref{d}) and (\ref{t}) yield
$$d_N -\tilde d_N = \sigma^2 Z_N + O((t_N-\tilde t_N)^2)
\mbox{~~and~~} t_N-\tilde t_N =O(Z_N),$$
where $$N^{2/3}Z_N=N^{2/3}\left\{g_{\mu_{A_N}}(\tilde t_N) -g_{ \frac{\tau}{N^{1/6}}}(\tilde t_N)\right\}$$
converges weakly to a centered  Gaussian distribution with variance $v^2/\sigma^4$.
Thus, it comes readily that 
\begin{equation}\label{dtilded} \sigma^{-1} N^{2/3} \left\{d_N-\tilde d_N\right\} \xrightarrow{\cal D}  N(0, v^2/\sigma^2).\end{equation}
Now (\ref{fluctuJohansson}) readily follows  since by Theorem \ref{theo: A}, 
given $A_N$, $\sigma^{-1}N^{2/3}\left(\lambda_{1}(M_N) -d_N\right)$ converges weakly to the Tracy-Widom distribution.

\item Another example is provided by 
\cite{LeeSchnelli} who considered the following deformed model 
$$\hspace*{-0.9cm}{W_N}+ {\rm diag}(v_1,\ldots,v_N)$$ where
$W_N$ is a Wigner matrix
and $v_i $ are {i.i.d}  random variables independent with $W_N$, with distribution $$ d\nu(x)= Z^{-1}(1 + x)^a(1 -x)^b f(x)1_{[-1,1]}(x)dx$$
 with $a<1,{b>1}$ and  $ f>0$ is a ${\cal C}^1$-function. 
Assume that $W_N$ is a G.U.E..
 Let 
 $\sigma_0$  be defined by {$\int \frac{1}{(1-x)^2}\mathrm{d}\nu(x)=\frac{1}{\sigma_0^2}$}. According to Section 3.2.1, we have  $ \mbox{ supp}\left(\mu_{sc}\boxplus \nu\right)=[d_\sigma^-,d_\sigma^+]$. Morevover, by (\ref{regular}), 
for all $\sigma> \sigma_0$, $d_\sigma^+ $ is a regular edge, ${ p(x)\sim C (d_\sigma^+ -x )^{\frac{1}{2}}}$.
Therefore denoting by $d_\sigma^+(N)$ the  (random) upper right edge of $ \mbox{ supp}\left(\mu_{sc}\boxplus \mu_{A_N}\right)$, Theorem \ref{theo: A} yields
$$\alpha^{-1}N^{2/3}(\lambda_{1}(M_N) -d_\sigma^+(N)) \xrightarrow{\cal D}TW.$$ A similar study as in Johansson's  example shows  that  the random edge  $d_\sigma^+(N)$ fluctuates 
as $$ \sqrt{N}(d_\sigma^+(N) -d_\sigma^+) \xrightarrow{\cal D} {\cal N}\left(0, \sigma^2\left(1-\sigma^2 \left(g_{\mu_{sc}\boxplus\nu}(d_\sigma^+)\right)^2\right)\right).$$ Thus,  we can deduce that $$ \sqrt{N}(\lambda_{1}(M_N) -d_\sigma^+) \xrightarrow{\cal D}  {\cal N}\left(0, \sigma^2\left(1-\sigma^2 \left(g_{\mu_{sc}\boxplus\nu}(d_\sigma^+)\right)^2\right)\right).$$

Note that for all $\sigma< \sigma_0$, according to (\ref{notregular}),
 ${ p(x)\sim C (d_\sigma^+ -x )^{b}}$.  Lee and Schnelli also investigate the fluctuations at the non-regular edge $d_\sigma^+$
and establish that 
$$ {N^{\frac{1}{b+1}}(\lambda_{1}(M_N) -d_\sigma^+) \xrightarrow{\cal D}G_{b+1}(s) }$$
as $N$ goes to infinity, where
$G_{b+1}(s) =(1-\exp({(\frac{s}{c})}^{b+1}))\1_{[0;+\infty[}(s)$
({Weibull distribution } with parameters $ b + 1$ and $c=c(\nu,\sigma)$).
\end{itemize}
\subsection{The general i.i.d case}
We finish by  some remarks concerning non Gaussian frameworks. We do not detail the results  since they do not fall under the scope of free probability theory.  
\begin{itemize}
\item  Universality at soft edges.\\
Several recent works
proved the  universality of the Tracy-Widom fluctuations at soft edges  for  quite general deformed Wigner matrices or sample covariance matrices without outliers.
The methods pursue
a Green function comparison strategy \cite{BPZ, univLS,univLS2} or make use of anisotropic local laws \cite{KY2}.
\item Non universal fluctuations  of outliers\\
 A new phenomenon arises  for the fluctuations of the outliers : the limiting distribution can depend on the distribution of the entries (non universality), according to the localization/delocalization of the eigenvectors of $A_N$. Note that in the Gaussian case in the previous subsection, the eigenvectors of the perturbation are irrelevant for the fluctuations, due to the unitary invariance in Gaussian models. 
We illustrate this dependence on the eigenvectors on $A_N$ in a very simple situation, in the additive case.
Consider two finite rank perturbations of rank 1, with one non null eigenvalue $\theta > \sigma$. The first one $A_N^{(1)}$ is a matrix with all entries equal to $\theta/N$ (delocalized eigenvector associated to $\theta$). The second one $A_N^{(2)}$ is a diagonal matrix (localized eigenvector). 
The fluctuations of the largest eigenvalue  $\lambda_1$ of the matrix $M_N^{(i)} = X_N + A_N^{(i)}$ ($i=1,2$) around 
$\rho_{\theta} := \theta + \frac{\sigma^2}{\theta}$  are given as follows :
\begin{proposition} Fluctuations of outliers
\begin{enumerate} 
\item Delocalized case \cite{FP} : The largest eigenvalue $\lambda_1(M^{(1)}_N)$ have Gaussian fluctuations :
\begin{equation} \label{delocalized}
\sqrt{N}  (\lambda_1(M_N^{(1)}) - \rho_\theta) \stackrel{{\cal D}}{\longrightarrow}  {\cal N}(0, \sigma^2( 1 - \sigma^2/ \theta^2))
\end{equation}
\item Localized case \cite{CDF} :  The largest eigenvalue $\lambda_1(M^{(2)}_N)$ fluctuates as 
\begin{equation} \label{localized}
 \sqrt{N} (1 - \frac{\sigma^2}{\theta^2}) (\lambda_1(M_N^{(2)}) - \rho_\theta) \stackrel{{\cal D}}{\longrightarrow} \mu \star {\cal N}(0, v_\theta).
 \end{equation}
where $\mu$ is the distribution of the entries of the Wigner matrix, the variance $v_\theta$ of the Gaussian distribution depends on $\theta$ and the second and fourth moments of $\mu$.
\end{enumerate}
\end{proposition}
The proof of \eqref{delocalized} is combinatorial and is based on the computation of large moments of $\Tr(M_N)$. 
The proof of \eqref{localized} relies on a determinant identity, analogous to \eqref{Syl}, boiling down to the behavior of a fixed rank determinant and a CLT for quadratic forms.

We refer to \cite{BY} (sample covariance case), \cite{CDF2, RenfrewSos1, RenfrewSos2, KY} for fluctuations of ouliers for more general perturbations, in the case of iid models, and \cite{BGM} for unitarily invariant models.\end{itemize}



\end{document}